\documentclass[12pt]{amsart}
\usepackage{amssymb}
\usepackage{amsfonts}
\usepackage{amsmath}
\usepackage{array}
\usepackage{rotating}
\usepackage{epsf}
\usepackage{tikz}
\usetikzlibrary{positioning,arrows,matrix,chains,shadows,shapes,calc,trees}
\tikzstyle{Point} = [fill, radius=0.08]
\tikzset{baseline={([yshift=-3.5pt]current bounding box.center)}}
\tikzset{level distance = 0.7cm, sibling distance = 2em}
\tikzset{edge from parent/.style={
       draw,
       edge from parent path = {(\tikzparentnode) --
                                (\tikzchildnode)}
    }}
\tikzset{root/.style={}}
\tikzset{white/.style={}}

\usepackage[matrix,frame,xdvi]{xypic} % ne pas enlever

\newdimen\vcadre\vcadre=0.1cm % marges verticales de la boite
\newdimen\hcadre\hcadre=0.1cm % marges horizontales de la boite
\def\GrTeXBox#1{\vbox{\vskip\vcadre\hbox{\hskip\hcadre%
      $#1$%
   \hskip\hcadre}\vskip\vcadre}}
\def\arx#1[#2]{\ifcase#1 \relax \or%
  \ar @{-}[#2]  \or%
  \ar @2{-}[#2] \or%
  \ar @{--}[#2] \or%
  \ar @2{.}[#2] \or%
  \ar @{~}[#2]  \fi}

\newtheorem{example}{Example}[section]

\newtheorem{note}[example]{Note}
\newtheorem{theorem}[example]{Theorem}

\newtheorem{definition}[example]{Definition}

\def\Proof{\noindent \it Proof -- \rm}
\def\qed{\hspace{3.5mm} \hfill \vbox{\hrule height 3pt depth 2 pt width 2mm}
\bigskip}

\def\sk{{\rm sk}}
\def\gr{{\rm gr}}

\def\evt{{\rm pEv}}
\def\E{{\mathbb E}}

\def\PQSym{{\bf PQSym}}
\def\CQSym{{\bf CQSym}}
\def\SQSym{{\bf SQSym}}
\def\ev{{\rm Ev}}

\def\park{{\bf a}}

\def\<{\langle}
\def\>{\rangle}
\def\R{{\mathbb R}}
\def\C{\operatorname{\mathbb C}}

\def\N{\operatorname{\mathbb N}}

\def\F{{\bf F}}
\def\S{{\bf S}}

\def\B{{\bf B}}

\def\SG{{\mathfrak S}}

\def\Sym{{\bf Sym}}
\def\NCSF{{\bf Sym}}

\def\maj{{\rm maj}}
\def\sep{\,|\,}   
\def\PS{{\bf P}}      % fonctions P de Schroder
\def\tg{\tilde g}

\def\PF{{\rm PF}}

\def\Tabvrule{\vrule width-0.4pt}       % Difference de largeur
\def\Tabhrule{\hrule \hrule height-0.4pt} % Difference de hauteur
\def\Tabstrut{\vrule height2.2ex % Sur la ligne
                     depth0.8ex  % Sous la ligne
                     width0ex    % centrage horizontal
\relax}

\def\PasCase#1{\omit%
            $\vcenter{\hbox {\vbox to 0.4pt{}}
               \hbox{\makebox[3ex]{\Tabstrut$#1$}}}%
               \Tabvrule$}
\def\PasCasePoint{\PasCase{\cdot}}
\def\DessinCarre#1{%
    \vcenter{\hbox{}\hrule
             \hbox{\vrule\makebox[3ex]{\Tabstrut$#1$}\vrule}\Tabhrule}%
             \Tabvrule}
\def\GenRuban#1{\vcenter{\halign{&$\DessinCarre{##}$\cr#1}}\egroup}

\def\sTabvrule{\vrule width-0.4pt}
\def\sTabhrule{\hrule \hrule height-0.4pt}
\def\sTabstrut{\vrule height1.6ex depth0.6ex width0ex \relax}
\def\sDessinCarre#1{%
    \vcenter{\hbox{}\hrule
             \hbox{\vrule\makebox[2.3ex]%
                  {\sTabstrut$\scriptstyle#1$}\vrule}\sTabhrule}%
             \sTabvrule}
\def\sGenRuban#1{\vcenter{\halign{&$\sDessinCarre{##}$\cr#1}}\egroup}

\def\ruban{%
  \bgroup
  \let\ =\omit
  \let\\=\cr
  \let\x=\times
  \let\.=\PasCasePoint
  \offinterlineskip
  \GenRuban}

\def\sruban{%
  \bgroup
  \let\ =\omit
  \let\x=\times
  \let\\=\cr
  \offinterlineskip
  \sGenRuban}
\def\Xint#1{\mathchoice
   {\XXint\displaystyle\textstyle{#1}}%
   {\XXint\textstyle\scriptstyle{#1}}%
   {\XXint\scriptstyle\scriptscriptstyle{#1}}%
   {\XXint\scriptscriptstyle\scriptscriptstyle{#1}}%
   \!\int}
\def\XXint#1#2#3{{\setbox0=\hbox{$#1{#2#3}{\int}$}
     \vcenter{\hbox{$#2#3$}}\kern-.5\wd0}}

\def\dashint{\Xint-}

\def\arbsu#1#2#3#4#5#6{
\xymatrix@R=0.1cm@C=2mm{
 & & & & {\GrTeXBox{#1}}\arx1[llld]\arx1[rrrd]\\
& {\GrTeXBox{#2}}\arx1[dl]\arx1[dr] &&&&&&
  {\GrTeXBox{#3}}\arx1[dl]\arx1[dr] \\
{\GrTeXBox{}} && {\GrTeXBox{#4}}\arx1[dl]\arx1[dr]
  &&&& {\GrTeXBox{#5}}\arx1[dl]\arx1[dr] && {\GrTeXBox{}}\\
& {\GrTeXBox{#6}}\arx1[dl]\arx1[dr] && {\GrTeXBox{}} && {\GrTeXBox{}}
&& {\GrTeXBox{}} \\
{\GrTeXBox{}} && {\GrTeXBox{}} \\
}
}

\def\arbsd#1#2#3#4#5{
\newdimen\vcadre\vcadre=0.01cm % marges verticales de la boite
\newdimen\hcadre\hcadre=0.01cm % marges horizontales de la boite
\xymatrix@R=0.1cm@C=2mm{
&&& {\GrTeXBox{#1}}\arx1[llld]\arx1[rrrd]\arx1[d]\\
{\GrTeXBox{}} &&&
  {\GrTeXBox{#2}}\arx1[dl]\arx1[dr] &&&
  {\GrTeXBox{#3}}\arx1[dl]\arx1[dr] \\
&& {\GrTeXBox{#4}}\arx1[dl]\arx1[dr]
  && {\GrTeXBox{}} & {\GrTeXBox{#5}}\arx1[dl]\arx1[dr] && {\GrTeXBox{}}\\
& {\GrTeXBox{}} && {\GrTeXBox{}} & {\GrTeXBox{}} && {\GrTeXBox{}} \\
}
}

\def\arbst#1#2#3#4#5{
\newdimen\vcadre\vcadre=0.01cm % marges verticales de la boite
\newdimen\hcadre\hcadre=0.01cm % marges horizontales de la boite
\xymatrix@R=0.1cm@C=2mm{
&&&& {\GrTeXBox{#1}}\arx1[llld]\arx1[rrrd] \\
& {\GrTeXBox{#2}}\arx1[dl]\arx1[dr] &&&&&&
  {\GrTeXBox{#3}}\arx1[dl]\arx1[dr] \\
{\GrTeXBox{}} && {\GrTeXBox{#4}}\arx1[dl]\arx1[dr]\arx1[d]
  &&&& {\GrTeXBox{#5}}\arx1[dl]\arx1[dr] && {\GrTeXBox{}} \\
& {\GrTeXBox{}} & {\GrTeXBox{}} & {\GrTeXBox{}} &&& {\GrTeXBox{}}
&& {\GrTeXBox{}} \\
}
}

\def\arbsq#1#2#3#4#5{
\newdimen\vcadre\vcadre=0.01cm % marges verticales de la boite
\newdimen\hcadre\hcadre=0.01cm % marges horizontales de la boite
\xymatrix@R=0.1cm@C=2mm{
 & & & & {\GrTeXBox{#1}}\arx1[llld]\arx1[rrrd]\\
& {\GrTeXBox{#2}}\arx1[dl]\arx1[dr] &&&&&&
  {\GrTeXBox{#3}}\arx1[dl]\arx1[d]\arx1[dr] \\
{\GrTeXBox{}} && {\GrTeXBox{#4}}\arx1[dl]\arx1[dr]
  &&&& {\GrTeXBox{}} & {\GrTeXBox{}} & {\GrTeXBox{}}\\
& {\GrTeXBox{#5}}\arx1[dl]\arx1[dr] &&\\
{\GrTeXBox{}} && {\GrTeXBox{}} \\
}
}

\def\arbsc#1#2#3#4{
\newdimen\vcadre\vcadre=0.01cm % marges verticales de la boite
\newdimen\hcadre\hcadre=0.01cm % marges horizontales de la boite
\xymatrix@R=0.1cm@C=2mm{
&&& {\GrTeXBox{#1}}\arx1[llld]\arx1[rrrd]\arx1[d]\\
{\GrTeXBox{}} &&&
  {\GrTeXBox{#2}}\arx1[dl]\arx1[dr]\arx1[d] &&&
  {\GrTeXBox{#3}}\arx1[dl]\arx1[dr] \\
&& {\GrTeXBox{}} & {\GrTeXBox{}} & {\GrTeXBox{}}
 & {\GrTeXBox{#4}}\arx1[dl]\arx1[dr] && {\GrTeXBox{}}\\
&&&& {\GrTeXBox{}} && {\GrTeXBox{}} \\
}
}

\def\arbss#1#2#3#4{
\newdimen\vcadre\vcadre=0.01cm % marges verticales de la boite
\newdimen\hcadre\hcadre=0.01cm % marges horizontales de la boite
\xymatrix@R=0.1cm@C=2mm{
&&& {\GrTeXBox{#1}}\arx1[llld]\arx1[rrrd]\arx1[d]\\
{\GrTeXBox{}} &&&
  {\GrTeXBox{#2}}\arx1[dl]\arx1[dr] &&&
  {\GrTeXBox{#3}}\arx1[dl]\arx1[d]\arx1[dr] \\
&& {\GrTeXBox{#4}}\arx1[dl]\arx1[dr]
  && {\GrTeXBox{}} & {\GrTeXBox{}} && {\GrTeXBox{}}\\
& {\GrTeXBox{}} && {\GrTeXBox{}} & {\GrTeXBox{}} && {\GrTeXBox{}} \\
}
}

\def\arbsse#1#2#3#4{
\newdimen\vcadre\vcadre=0.01cm % marges verticales de la boite
\newdimen\hcadre\hcadre=0.01cm % marges horizontales de la boite
\xymatrix@R=0.1cm@C=2mm{
&&&& {\GrTeXBox{#1}}\arx1[llld]\arx1[rrrd] \\
& {\GrTeXBox{#2}}\arx1[dl]\arx1[dr] &&&&&&
  {\GrTeXBox{#3}}\arx1[dl]\arx1[d]\arx1[dr] \\
{\GrTeXBox{}} && {\GrTeXBox{#4}}\arx1[dl]\arx1[dr]\arx1[d]
  &&&& {\GrTeXBox{}} && {\GrTeXBox{}} \\
& {\GrTeXBox{}} & {\GrTeXBox{}} & {\GrTeXBox{}} &&& {\GrTeXBox{}}
&& {\GrTeXBox{}} \\
}
}

\def\arbsh#1#2#3{
\newdimen\vcadre\vcadre=0.01cm % marges verticales de la boite
\newdimen\hcadre\hcadre=0.01cm % marges horizontales de la boite
\xymatrix@R=0.1cm@C=2mm{
&&& {\GrTeXBox{#1}}\arx1[llld]\arx1[rrrd]\arx1[d]\\
{\GrTeXBox{}} &&&
  {\GrTeXBox{#2}}\arx1[dl]\arx1[dr]\arx1[d] &&&
  {\GrTeXBox{#3}}\arx1[dl]\arx1[dr]\arx1[d] \\
&& {\GrTeXBox{}} & {\GrTeXBox{}} & {\GrTeXBox{}}
 & {\GrTeXBox{}} & {\GrTeXBox{}} & {\GrTeXBox{}}\\
}
}

\newenvironment{arb}{\begin{tikzpicture}[%baseline,
scale=0.5,level distance=7mm,level 1/.style={sibling distance=8mm},level 2/.style={sibling distance=5mm},level 3/.style={sibling distance=5mm},grow'=down, font=\scriptsize]
\tikzstyle{ve}=[draw,circle,inner sep=1pt,fill] 
\tikzstyle{vv}=[draw,circle,inner sep=1pt] 
\tikzstyle{vee}=[minimum size=0pt ,inner sep=0pt]}{\end{tikzpicture}}

\newcommand{\rd}[1]{\node[ve,label=above:$#1$] {}}
\newcommand{\vv}{node[vv] {}}
\newcommand{\ve}{node[ve] {}}

\newcommand{\vl}[1]{node[ve,label=right:$#1$] {}}
\newcommand{\vr}[1]{node[ve,label=left:$#1$] {}}

\newcommand{\va}[1]{node[ve,label=below:$#1$] {}}

\def\arbA{{\scriptstyle \circ}}

\def\arbB{\begin{arb}
\node[ve] {}
child{\vv} child{\vv};
\end{arb}}

\def\arbCA{\begin{arb}
\node[ve] {}
child{\ve  child{\vv} child{\vv}} child{\vv};
\end{arb}}

\def\arbCB{\begin{arb}
\node[ve] {}
child{\vv} child{\ve  child{\vv} child{\vv}};
\end{arb}}

\def\arbCC{\begin{arb}
\node[ve] {}
child{\vv} child{\vv} child{\vv};
\end{arb}}

\def\arbD{\begin{arb}
\node[ve] {}
child{\vv} child{\vv} child{\vv} child{\vv};
\end{arb}}

\def\arbF{\begin{arb}
\node[ve] {}
child{\ve  child{\vv} child{\vv}} child{\ve  child{\vv} child{\ve  child{\vv} child{\vv} child{\vv} child{\vv}}} child{\vv};
\end{arb}}

\def\arbE{\begin{arb}
\node[ve] {}
child{\vv} child{\ve  child{\vv} child{\vv}} child{\vv};
\end{arb}}

\def\carbA{{\scriptstyle \bullet}}

\def\carbUn{\begin{arb}
\node[ve] {}
child{\ve};
\end{arb}}

\def\carbBA{\begin{arb}
\node[ve] {}
child{\ve} child{\ve};
\end{arb}}

\def\carbBB{\begin{arb}
\node[ve] {}
child{\ve child{\ve}};
\end{arb}}

\def\carbCA{\begin{arb}
\node[ve] {}
child{\ve} child{\ve} child{\ve};
\end{arb}}

\def\carbCB{\begin{arb}
\node[ve] {}
child{\ve} child{\ve  child{\ve}};
\end{arb}}

\def\carbCC{\begin{arb}
\node[ve] {}
child{\ve child{\ve}} child{\ve};
\end{arb}}

\def\carbCD{\begin{arb}
\node[ve] {}
child{\ve child{\ve} child{\ve}};
\end{arb}}

\def\carbF{\begin{arb}
\node[ve] {}
child{\ve  child{\ve} child{\ve}}  child{\ve} child{\ve  child{\ve}} ;
\end{arb}}

\def\carbCE{\begin{arb}
\node[ve] {}
child{\ve child{\ve  child{\ve }}};
\end{arb}}

\title[Combinatorics of Poincar\'e's and Schr\"oder's equations]%
{Combinatorics of Poincar\'e's \\and Schr\"oder's equations}

\author[F. Menous, J.-C.~Novelli and J.-Y.~Thibon]%
{Fr\'ed\'eric Menous, Jean-Christophe Novelli and Jean-Yves Thibon}
\address[Menous]{Laboratoire de Math\'ematiques, B\^atiment 425, Universit\'e
Paris-Sud, 91405 Orsay Cedex, FRANCE}
\email[Fr\'ed\'eric Menous]{frederic.menous@math.u-psud.fr}
\address[Novelli and Thibon] {Universit\'e Paris-Est Marne-la-Vall\'ee \\
Laboratoire d'Infor\-ma\-tique Gaspard-Monge (CNRS - UMR 8049)\\
77454 Marne-la-Vall\'ee cedex 2 \\
FRANCE}
\email[Jean-Christophe Novelli]{novelli@univ-mlv.fr}
\email[Jean-Yves Thibon]{jyt@univ-mlv.fr} 


\begin{document}

\begin{abstract}
We investigate the combinatorial properties of the functional
equation $\phi[h(z)]=h(qz)$ for the conjugation of a formal
diffeomorphism $\phi$ of $\C$ to its linear part $z\mapsto qz$.
This is done by interpreting the functional equation in terms of
symmetric functions, and then lifting it to noncommutative symmetric functions.
We describe explicitly the expansion of the solution in terms of
plane trees and prove that its expression on the ribbon basis has
coefficients in $\N[q]$ after clearing the denominators $(q)_n$.
We show that the conjugacy equation can be lifted to a quadratic
fixed point equation in the free triduplicial algebra on one generator.
This can be regarded as a $q$-deformation of the duplicial interpretation
of the noncommutative Lagrange inversion formula. Finally, these calculations
are interpreted in terms of the group of the operad of Stasheff polytopes, and
are related to Ecalle's arborified expansion by means of morphisms between
various Hopf algebras of trees.
\end{abstract}

\maketitle
\tableofcontents
%%%%%%%%%%%%%%%%%%%%%%%%%%%%%%%%%%%%%%%%%%%%%%%%%%%%%%%%%%%%%%%%%%%%%%%%%%%%%%%
%%%%%%%%%%%%%%%%%%%%%%%%%%%%%%%%%%%%%%%%%%%%%%%%%%%%%%%%%%%%%%%%%%%%%%%%%%%%%%%
%%%%%%%%%%%%%%%%%%%%%%%%%%%%%%%%%%%%%%%%%%%%%%%%%%%%%%%%%%%%%%%%%%%%%%%%%%%%%%%
\section{Introduction}

Algebraic identities between generic formal power series can often
be interpreted as identities between symmetric functions.
This is the case, for example, with the Lagrange inversion
formula 
(see, e.g., \cite{Mcd}, Ex. 24 p. 35, Ex. 25 p. 132, \cite{Las}
Section~2.4, and \cite{Len}). The problem can be stated as follows.
Given 
\begin{equation}
\varphi(z)=\sum_{n\ge 0}\varphi_nz^n\quad (\varphi_0\not=0)\,
\end{equation}
find the coefficients $g_n$ of the unique power series 
\begin{equation}
g(z)=\sum_{n\geq0}g_nz^{n+1}
\quad\text{satisfying}\quad
z = \frac{g(z)}{\varphi(g(z))} \,.
\end{equation}
We can assume that $\varphi_0=1$ and that
\begin{equation}
\varphi(g)=\sum_{n\ge 0}h_n(X)g^n= \prod_{n\ge 1}(1-gx_n)^{-1}=: \sigma_g(X)
\end{equation}
is the generating series of the homogeneous symmetric functions of an
infinite set of variables $X$. 
In  $\lambda$-ring notation, the solution reads
\begin{equation}
g_n =\frac1{n+1}h_n((n+1)X)
\end{equation}
(recall that $\sigma_t(nX)=\sigma_t(X)^n$, see, e.g., \cite{Mcd} p. 25).
On this expression, it is clear that $g_n$ is Schur positive, in fact, it is
the Frobenius characteristic of the  permutation representation of $\SG_n$
on the set $\PF_n$ of  parking functions of length $n$. 
These calculations can be lifted to the algebra of noncommutative symmetric
functions, and the result is then interpreted in terms of representations of
0-Hecke algebras. This in turn leads to various combinatorial interpretations,
to $q$-analogues, and to a new interpretation of the antipode of the Hopf
algebra of noncommutative formal diffeomorphisms \cite{NTLag}.

There is another functional equation which can be investigated in this setting.
Given a formal diffeomorphism
\begin{equation}
\phi(z) = \sum_{n\ge 0}\phi_n z^{n+1}\quad \text{with $\phi_0=q\not=0$},
\end{equation}
one may look for a formal diffeomorphism 
tangent to identity 
\begin{equation}
h(z)=\sum_{n\ge 0}g_nz^{n+1} = zg(z)\quad (g_0=1),
\end{equation}
conjugating $\phi$ to its linear part
\begin{equation}\label{eq:conjug}
h^{-1}\circ\phi\circ h(z)=qz\ \text{or equivalently}
                  \  \phi[h(z)]=h(qz)=qzg(qz).
\end{equation}
In terms of symmetric functions, we can assume that
\begin{equation}
\phi(z) = qz\sigma_z(X)
\end{equation}
so that the conjugacy equation reads
\begin{equation}
\phi[h(z)] = qh(z)\sigma_{h(z)}(X) = qz\sum_{n\ge 0}g_n (qz)^n,
\end{equation}
and interpreting $g_n$ as symmetric functions $g_n(X)$, we can get rid of $z$
by homogeneity (since $g_n(zX)=z^ng_n(X)$). Our functional equation reads now
\begin{equation}
g(X)\sigma_{g(X)} = g(qX).
\end{equation}
We can lift this to noncommutative symmetric functions, for example as
\begin{equation}\label{eq:ncconj}
g(qA) = \sum_{n\ge 0} S_n(A) g(A)^{n+1}.
\end{equation}
For $q=0$, this reduces to the functional equation for the antipode of
the noncommutative Fa\`a di Bruno Hopf algebra \cite{BFK,NTLag}, so that this problem can indeed be
regarded as a generalisation of the noncommutative Lagrange inversion.

The conjugacy equation for $h$ is often called Poincar\'e's equation, and the
equivalent one for $h^{-1}$, Schr\"oder's equation. Indeed, it has been first
discussed by Schr\"oder \cite{Schr2}, who discovered a few explicit solutions,
which are still essentially the only known ones. It is easy to show the existence and
unicity of  a formal solution when $q$ is not a root of unity. The
analyticity of the solution for $|q|\not=1$ has been established by Koenigs
\cite{Koe}. It is interesting that this result can be easily proved by means
of inequalities involving  the Schr\"oder numbers \cite{Mar}, defined by the
same Schr\"oder in a totally different context \cite{Schr1}. Much more
difficult is Siegel's proof of convergence in the case $q=e^{2\pi i\theta}$
with $\theta$ satisfying a diophantine condition \cite{Sie}
(see also \cite{smf} for a modern proof under Bruno's condition). Again in this
case, the Schr\"oder numbers play a crucial role in the majorations.

We shall see that analyzing the conjugacy equation at the level of the
noncommutative Fa\`a di Bruno algebra provides a simple explanation of this
fact, by letting Schr\"oder trees appear naturally in the iterative
solution of a $q$-difference equation. 
The resulting expressions turn out to be identical to those produced by Ecalle's
arborification method \cite{JESNAG,FM,smf}. This coincidence will be explained in
Section \ref{cvsnc}, where it will be proved that both methods can be interpreted
in terms of  calculations in the group of an  operad and in related Hopf algebras.

Identifying the noncommutative Fa\`a di Bruno algebra with noncommutative
symmetric functions as in \cite{NTLag}, we have several bases at our disposal. The solution $g$ of
the noncommutative Poincar\'e equation is naturally expressed in the complete
basis $S^I$. After clearing out the denominators $(q;q)_n$, it turns out that
its homogeneous components $g_n$ are positive on the ribbon basis. This unexpected
fact suggests that these should be the graded characteristics of some projective
modules over $0$-Hecke algebra, a conjecture that we expect to investigate in another
paper. This positivity property will be proved in two different ways.
We shall first recast the conjugacy equation as a quadratic fixed point problem,
by means of the triduplicial operations introduced in \cite{NTDup}. On the ribbon
basis, the quadratic map is manifestly positive. Next, comparing the binary tree expansion
with the previous one based on reduced plane trees, we obtain a natural bijection
between  these trees and hypoplactic classes of parking functions ({\it aka} parking quasi-ribbons
or segmented nondecreasing parking functions). This solves a problem which was left
open in \cite{NT1}, and provides a bijection similar to the duplicial bijection
of \cite{NTDup} between nondecreasing parking functions and binary trees.

In Section \ref{sec:rib}, we describe the expansion of $g_n$
on the ribbon basis. The numerator of each coefficient is a $q$-analogue
of $n!$, recording a statistic on permutations which is explicitly described.

In Section \ref{schrodeq}, we discuss Schr\"oder's equation at the level  of
noncommutative symmetric functions. It leads to
a different combinatorics. There is no natural expansion on trees, but instead,
there is a rather explicit algebraic formula for the coefficients, which amounts
to applying a simple linear transformation to a famous sequence of noncommutative
symmetric functions, the $q$-Klyachko elements $K_n(q)$ \cite{NCSF1,NCSF2}, which
occur as well as Lie idempotents in descent algebras or as noncommutative Hall-Littlewood
functions \cite{Hiv}. It is then shown in Section \ref{ncmould} that the same coefficients arise
when the problem is considered from the point of view of mould calculus and differential operators. 
Thus, at least for this problem, the mould calculus approach can be seen to be  dual to that
relying on the noncommutative Fa\`a di Bruno algebra.

The rest of the paper is devoted to the explanation of the coincidence between
the coefficients of our first plane tree expansion, and Ecalle's arborified
coefficients. The short story is that on the one hand, the paper \cite{FM} provides
an interpretation of the arborification method as a lift of the original problem
to an equation in the group of characters of a Connes-Kreimer algebra. On the other
hand, our version with plane trees  of the functional equation can be naturally
interpreted in the group of a free operad. This group turns out to be isomorphic
to the group of characters of a Hopf algebra of reduced plane forests, which 
admits a surjective morphism to the previous Connes-Kreimer algebra.

Section \ref{opred} provides some background on the operad of reduced plane trees.
It is a free operad with one generator in each degree $n\ge 2$, also known as
the operad of Stasheff polytopes, or as a free $S$-magmatic operad \cite{RH,JLL}.
We describe the associated group, and prove that it is isomorphic to the group
of characters of the Hopf algebra of reduced plane trees of \cite{MEF}.

In Section \ref{ncfd}, we explain the encoding of
the previous group by means of Polish codes of trees, and illustrate the method
on the cases of Lagrange inversion and of the Poincar\'e equation.

In Section \ref{cvsnc}, we recall the Hopf algebraic interpretation
of the arborification method \cite{FM,smf}, and prove that the skeleton map
already introduced in \cite{NTLag} induces a morphism of Hopf algebras
between reduced plane forests and the $\N^*$-decorated Connes-Kreimer algebra.

In Section \ref{sec:cat}, we review briefly the interpretation
of Lagrange inversion and of Cayley's formula for the solution of a generic
differential equation in terms of an operad on (non-reduced) plane trees.

%%%%%%%%%%%%%%% 
Finally, it is generally interesting to look at the images of  formal
series in combinatorial Hopf algebras under various characters.
In the Appendix (Section \ref{sec:ex}), we review a few examples of explicit
solutions of the conjugacy equation.
Apart from the trivial case of linear fractional transformations
$\phi(z)=qz/(1-z)$ (corresponding to the alphabet $A=\{1\}$), there is the
already nontrivial case of the logistic map $\phi(z)=qz(1-z)$, 
corresponding to $A=\{-1\}$, for which explicit solutions (already given by
Schr\"oder) are known for $q=-2,2,4$.
The case $A={\mathbb E}$, corresponding to $\phi(z)=qze^z$ 
is not explicitly solved, but it leads to interesting statistics on pairs of
permutations. These examples are investigated numerically in \cite{CZ1,CZ2,CJZ,CV}.

\medskip
{\footnotesize
{\it Acknowledgements.-}
This project has been partially supported by
the project CARMA of the French Agence Nationale
de la Recherche.
}

%%%%%%%%%%%%%%%%%%%%%%%%%%%%%%%%%%%%%%%%%%%%%%%%%%%%%%%%%%%%%%%%%%%%%%%%%%%%%%%
%%%%%%%%%%%%%%%%%%%%%%%%%%%%%%%%%%%%%%%%%%%%%%%%%%%%%%%%%%%%%%%%%%%%%%%%%%%%%%%
%%%%%%%%%%%%%%%%%%%%%%%%%%%%%%%%%%%%%%%%%%%%%%%%%%%%%%%%%%%%%%%%%%%%%%%%%%%%%%%
\section{Notations}

This paper is a continuation of \cite{NTLag,NTDup}. 
Our notations for ordinary symmetric functions are as in \cite{Mcd}, and
for noncommutative symmetric functions as in \cite{NCSF1,NCSF2}. 

The classical algebra of symmetric functions, denoted by $Sym$ or $Sym(X)$, is
a free associative and commutative graded algebra with one generator in each degree:
\begin{equation}
Sym = \C[h_1,h_2,\ldots] = \C[e_1,e_2,\ldots] = \C[p_1,p_2,\ldots]
\end{equation}
where the $h_n$ are the complete homogeneous symmetric functions, the $e_n$ the elementary
symmetric functions, and the $p_n$ the power sums. 

Its usual bialgebra structure is defined by the coproduct
\begin{equation}
\Delta_0 h_n = \sum_{i=0}^nh_i\otimes h_{n-i}\quad (h_0=1)
\end{equation}
which allows to interpret it as the algebra of polynomial functions on the multiplicative group
\begin{equation}
G_0 = \{a(z)=\sum_{n\ge 0}a_n z^n\ (a_0=1)\}
\end{equation}
of formal power series with constant term 1: $h_n$ is the coordinate function 
\begin{equation}
h_n:\ a(z)\longmapsto a_n.
\end{equation}
Indeed, $h_n(a(z)b(z))=(\Delta_0 h_n) (a(z)\otimes b(z))$.

But $h_n$ can also be interpreted as a coordinate on the group
\begin{equation}
G_1 = \{A(z)=\sum_{n\ge 0}a_n z^{n+1}\ (a_0=1)\}
\end{equation}
of formal diffeomorphisms tangent to identity, under functional composition.
Again with $h_n(A(z))=a_n$ and $h_n(A(z)B(z))=\Delta_1(A(z)\otimes B(z))$,
the coproduct is now
\begin{equation}
\Delta_1 h_n = \sum_{i=0}^nh_i\otimes h_{n-i}((i+1)X)\quad (h_0=1)
\end{equation}
where $h_n(mX)$ is defined as the coefficient of $t^n$ in $(\sum h_kt^k)^m$.
The resulting bialgebra is known as the Fa\`a di Bruno algebra \cite{JR}.

These constructions can be repeated word for word with the algebra
$\Sym$ of noncommutative symmetric functions. It is a free associative
(and noncommutative) graded algebra with one generator $S_n$ in each degree,
which can be interpreted as above if the coefficients $a_n$ belong to
a noncommutative algebra. In this case, $G_0$ is still a group, but $G_1$
is not, as its composition is not anymore associative. However, the coproduct $\Delta_1$
\begin{equation}
\Delta_1 S_n = \sum_{i=0}^nS_i\otimes S_{n-i}((i+1)A)\quad (S_0=1)
\end{equation}
remains coassociative, and $\Sym$ endowed with this coproduct is a Hopf algebra,
known as Noncommutative Formal Diffeomorphims \cite{BFK,NTLag}, or as the noncommutative
Fa\`a di Bruno algebra \cite{EFLM}.

The classical trick of regarding a generic series as a series of symmetric functions
amounts to working in one of these Hopf algebras. The occurence of trees in
the solutions of certain problems can be traced back to the existence of Hopf
algebras morphisms between these algebras and various Hopf algebras of trees.

Recall that bases of $\Sym_n$ are labelled by
compositions $I$ of $n$. The noncommutative complete and elementary functions
are denoted by $S_n$ and $\Lambda_n$, and the notation $S^I$ means
$S_{i_1}\cdots S_{i_r}$. The ribbon basis is denoted by $R_I$.
The notation $I\vDash n$ means that $I$ is a composition of $n$.
The conjugate composition is denoted by $I^\sim$. The product formula for
ribbons is
\begin{equation}
R_IR_J = R_{IJ}+R_{I\triangleright J}
\end{equation}
where for $I=(i_1,\ldots,i_r)$ and $J=(j_1,\ldots,j_s)$,
\begin{equation}
IJ=(i_1,,\ldots,i_r,j_1,\ldots,j_s)\ \text{and}\ I\triangleright J=(i_1,,\ldots,i_r+j_1,\ldots,j_s).
\end{equation}

The graded dual of $\Sym$ is $QSym$ (quasi-symmetric functions).
The dual basis of $(S^I)$ is $(M_I)$ (monomial), and that of $(R_I)$
is $(F_I)$.

The \emph{evaluation} $\ev(w)$ of a word $w$ over a totally ordered alphabet
$A$ is the sequence $(|w|_a)_{a\in A}$ where $|w|_a$ is the number of
occurrences of $a$ in $w$. The {\em packed evaluation} $I=\evt(w)$ is the
composition obtained by removing the zeros in $\ev(w)$.

Two permutations $\sigma,\tau\in\SG_n$ are said to be sylvester-equivalent
if the decreasing binary trees\footnote{The decreasing tree $T(w)$
of a word without repeated letters $w=unv$ and maximal letter $n$
is the binary tree with root $n$ and left and right subtrees $T(u)$ and $T(v)$.}  
of $\sigma^{-1}$ and $\tau^{-1}$ have the same shape.
The generating function of the number of inversions on a sylvester class is given
by the $q$-hook-length formula \cite{BW,HNT}.

%%%%%%%%%%%%%%%%%%%%%%%%%%%%%%%%%%%%%%%%%%%%%%%%%%%%%%%%%%%%%%%%%%%%%%%%%%%%%%%
%%%%%%%%%%%%%%%%%%%%%%%%%%%%%%%%%%%%%%%%%%%%%%%%%%%%%%%%%%%%%%%%%%%%%%%%%%%%%%%
%%%%%%%%%%%%%%%%%%%%%%%%%%%%%%%%%%%%%%%%%%%%%%%%%%%%%%%%%%%%%%%%%%%%%%%%%%%%%%%
\section{Recursive solution of Poincar\'e's equation}

Equation~\eqref{eq:ncconj} can be written as a $q$-difference equation
\begin{equation}
g(qA)-g(A) = \sum_{n\ge 1}S_n(A)g(A)^{n+1}.
\end{equation}
Introducing a homogeneity parameter $t$, we have
\begin{equation}
\label{eqqt}
g(qtA)-g(tA) = \sum_{n\ge 1}t^n S_n(A)g(tA)^{n+1}.
\end{equation}
Let $g_n$ be the term of degree $n$ in $g$, so that
\begin{equation}
g(tA) = \sum_{n\ge 1}t^n g_n(A).
\end{equation}
Comparing the homogeneous components in
both sides of~\eqref{eqqt}, one gets a triangular system allowing to compute
the $g_n$ recursively. For $n\leq3$:
\begin{equation}
\begin{split}
g_0 &= 1, \\
(q-1) g_1 &= S_1, \\
(q^2-1) g_2 &= 2S_1g_1 + S_2, \\
(q^3-1) g_3 &= 2S_1g_2 + S_1g_1^2 + 3S_2g_1 + S_3.
\end{split}
\end{equation}
Define
\begin{equation}
q_n = q^n-1,\  (q)_n =q_n q_{n-1}\cdots q_1 \ \text{and}\ \tg_n = (q)_n g_n
\end{equation}
The first $\tg_n$ are then
\begin{equation}
\label{first-tg}
\begin{split}
\tg_1 &= S_1, \\
\tg_2 &= (q-1) S_2 + 2 S^{11}, \\
\tg_3 &= (q)_2 S_3 + 3(q^2-1)S^{21} + 2(q-1)S^{12} + (5+q)S^{111}.
\end{split}
\end{equation}
On the ribbon basis of $\NCSF$, the expression is quite remarkable:
\begin{equation}
\begin{split}
\tg_1 &= S_1, \\
\tg_2 &= (q+1) R_2 + 2 R_{11}, \\
\tg_3 &= (1+q)(1+q+q^2)R_3+ (2+q+3q^2)R_{21}+3(1+q)R_{12}+(5+q)R_{111}.
\end{split}
\end{equation}
Indeed, $\tg_n$ is a linear combination of ribbons with positive coefficients
which are all $q$-analogues of $n!$. Note that it is immediate, by induction on $n$, that $\tg_n|_{q=1}=n!S_1^n$, but
it is not clear that the coefficients are in $\N[q]$. A combinatorial interpretation of these
coefficients will be given below (Theorem~\ref{th:rib}).

%%%%%%%%%%%%%%%%%%%%%%%%%%%%%%%%%%%%%%%%%%%%%%%%%%%%%%%%%%%%%%%%%%%%%%%%%%%%%%%
%%%%%%%%%%%%%%%%%%%%%%%%%%%%%%%%%%%%%%%%%%%%%%%%%%%%%%%%%%%%%%%%%%%%%%%%%%%%%%%
%%%%%%%%%%%%%%%%%%%%%%%%%%%%%%%%%%%%%%%%%%%%%%%%%%%%%%%%%%%%%%%%%%%%%%%%%%%%%%%
\section{A tree-expanded solution}\label{sec:TreeExp}

In order to solve \eqref{eqqt}, 
define a $q$-integral by\footnote{This is just the ordinary $q$-integral up to conjugation
by the transformation $t\mapsto (q-1)t$.}
\begin{equation}
\dashint_a^b t^{n-1} d_qt =\left[\frac{t^n}{q^n-1}\right]_a^b 
\end{equation}
and a $q$-difference operator
\begin{equation}
\Delta_q f(t) =\frac{f(qt)-f(t)}{t}.
\end{equation}
Then,
\begin{equation}
\Delta_q \dashint_0^tf(s)d_q s = f(t)
\end{equation}
so that $g$ is the unique solution of the fixed point equation
\begin{equation}
g(tA) = g_0 + \sum_{n\ge 1}S_n(A)\dashint_0^t s^{n-1}g(sA)^{n+1}d_q s.
\end{equation}
This equation is of the form 
\begin{equation}
g = g_0 + \sum_{n\ge 2} F_n(g,\ldots,g)
\end{equation}
where 
\begin{equation}\label{eq:Fn}
F_n(x_1,\ldots,x_n) = S_{n-1}\dashint_0^1 s^{n-2}x_1(s)\cdots x_n(s) d_q s
\end{equation}
is an $n$-linear operator. 
The solution can therefore be expanded as a sum over
reduced plane trees (plane trees in which all internal vertices have at least
two descendants), which will be called  {\it Schr\"oder trees} in the sequel.

Proceeding as in \cite{NTLag}, we introduce another indeterminate $S_0$
(noncommuting with the other $S_n$) and set $g_0=S_0$. The solution
is then a linear combinations of monomials $S^I$ where $I$ is 
a vector of nonnegative integers, with $i_1>0$. 

The first $\tg_n$ are then
\begin{equation}
\label{first-tg0}
\begin{split}
\tg_1 &= S^{100}\\
\tg_2 &= (q-1)S^{2000} + S^{11000} + S^{10100} \\
\tg_3 &= (q)_2S^{30000}+(q^2-1)(S^{210000}+S^{201000}+S^{200100})\\
  &\quad+ (q-1)(S^{120000}+S^{102000})\\
  &\quad+ S^{1110000}+S^{1101000}+S^{1011000}+S^{1010100} + (q+1)S^{1100100}.
\end{split}
\end{equation}
We can interpret each $S_i$ as the symbol of an $(i+1)$-ary operation in
Polish notation.  Then,  $\tg_n$ is a sum over Polish codes of 
Schr\"oder trees as in \cite[Fig. 4]{NTLag}:

%%%%%%%%%%%%%%%%%%%%%%%%%%% A
\def\ftu{
%\newdimen\vcadre\vcadre=0.01cm % marges verticales de la boite
%\newdimen\hcadre\hcadre=0.01cm % marges horizontales de la boite
\xymatrix@R=0.1cm@C=2mm{
& {\GrTeXBox{2}}\arx1[ld]\arx1[rd]\arx1[d]\\
% ligne 2
{\GrTeXBox{0}} & {\GrTeXBox{0}} & {\GrTeXBox{0}} \\
}
}

\def\ftd{
%\newdimen\vcadre\vcadre=0.01cm % marges verticales de la boite
%\newdimen\hcadre\hcadre=0.01cm % marges horizontales de la boite
\xymatrix@R=0.1cm@C=2mm{
&& {\GrTeXBox{1}}\arx1[ld]\arx1[rd] \\
% ligne 2
& {\GrTeXBox{1}}\arx1[ld]\arx1[rd] && {\GrTeXBox{0}} \\
% ligne 3
{\GrTeXBox{0}} && {\GrTeXBox{0}} \\
}
}

\def\ftt{
%\newdimen\vcadre\vcadre=0.01cm % marges verticales de la boite
%\newdimen\hcadre\hcadre=0.01cm % marges horizontales de la boite
\xymatrix@R=0.1cm@C=2mm{
& {\GrTeXBox{1}}\arx1[ld]\arx1[rd] \\
% ligne 2
{\GrTeXBox{0}} && {\GrTeXBox{1}}\arx1[ld]\arx1[rd] \\
% ligne 3
& {\GrTeXBox{0}} && {\GrTeXBox{0}} \\
}
}

\begin{align}
\tg_2 &= (q-1)\ftu &+& \ftd &+& \ftt \\
    &= (q-1)S^{2000} &+&\qquad S^{11000}
       &+&\qquad S^{10100}.
\end{align}

The exponent vectors $I$ encoding Schr\"oder trees as above will be referred to
as {\it Schr\"oder pseudocompositions}.

From \eqref{eq:Fn}, we have: 

\begin{theorem}
Let $I$ be a Schr\"oder pseudocomposition, and $T(I)$ be the tree encoded by $I$. The coefficient of
$S^I$ in $g$ is
\begin{equation}
m_I(q) = \prod_{v\in T(I)}\frac1{q_{\phi(v)-1}}
\end{equation}
where $v$ runs over the internal vertices of $T(I)$ and $\phi(v)$ is the
number of leaves of the subtree of $v$.\qed
\end{theorem}

\begin{note}
{\rm These coefficients are precisely those obtained by Ecalle's
arborification method \cite{JESNAG,FM,smf}.
This coincidence will be explained in Section \ref{cvsnc}.
}
\end{note}
%%%%%%%%%%%%%%%%%%%%%%%%%%%%%%%%%% B
\begin{example}{\rm
For 
\def\bb#1#2#3#4#5{
%\newdimen\vcadre\vcadre=0.01cm % marges verticales de la boite
%\newdimen\hcadre\hcadre=0.01cm % marges horizontales de la boite
\xymatrix@R=0.1cm@C=2mm{
 & & & & {\GrTeXBox{#1}}\arx1[llld]\arx1[rrrd]\arx1[d]\\
% ligne 2
& {\GrTeXBox{#2}}\arx1[dl]\arx1[dr] &&& {\GrTeXBox{#3}}\arx1[dl]\arx1[dr]
&&& {\GrTeXBox{#5}} \\
% ligne 3
{\GrTeXBox{#5}} && {\GrTeXBox{#5}}
  & {\GrTeXBox{#5}}
 && {\GrTeXBox{#4}}\arx1[dl]\arx1[d]\arx1[dr]\arx1[drr] \\
% ligne 4
&&&& {\GrTeXBox{#5}} & {\GrTeXBox{#5}} & {\GrTeXBox{#5}} & {\GrTeXBox{#5}} \\
}
}

\begin{equation}
T= \bb{2}{1}{1}{3}{0}
\end{equation}
decorating each internal vertex $v$ with the number $\phi(v)-1$, we obtain 
\begin{equation}
\bb{7}{1}{4}{3}{}
\end{equation}
so that
\begin{equation}
m_{210010300000} = \frac{1}{q_7\,q_4q\,_3\,q_1}.
\end{equation}
}
\end{example}

\begin{note}\label{note:bin}{\rm
The $I$ whose nonzero entries are all equal to 1 correspond to binary trees.
The anti-refinements $J\preceq I$ of such an $I$, obtained by summing
consecutive nonzero entries in all possible ways, correspond to the trees
$T(J)$ obtained by contracting (internal) left edges in all possible ways in
$T(I)$. This procedure provides a way to group Schr\"oder trees into
classes labelled by binary trees. An algebraic interpretation of these groups
will be provided below.
}
\end{note}

\begin{example}\label{ex:contracts}{\rm
The  contractions of the binary tree $T(1101100011000)$ are

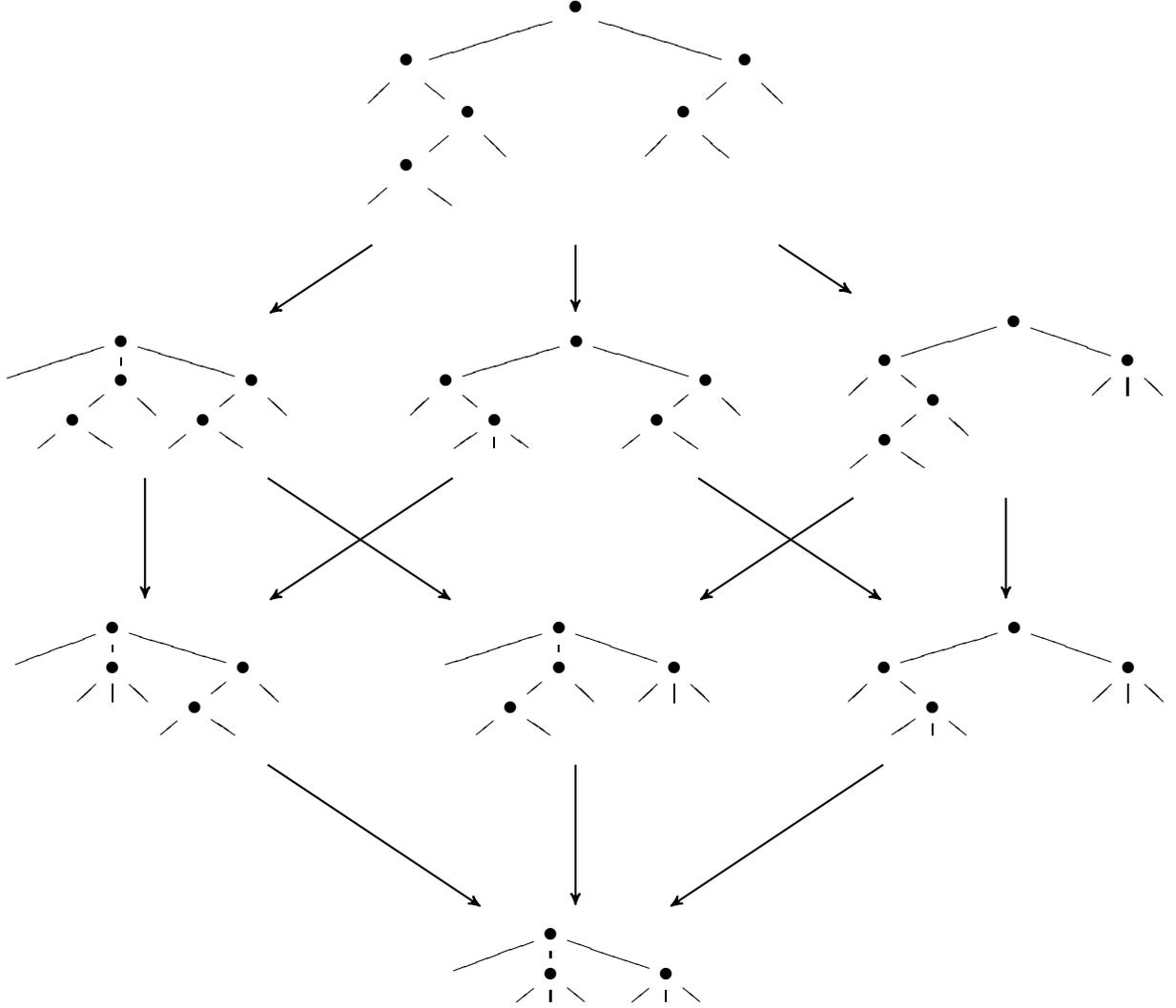
\begin{tikzpicture}[scale=.4,>=stealth',shorten >=1pt,auto,node distance=2cm,
                    thick,main node/.style={}]

\path (0,0) node (1)
    {\arbsu{\bullet}{\bullet}{\bullet}{\bullet}{\bullet}{\bullet}};
\path (-15,-10)  node (2)  %
    {\arbsd{\bullet}{\bullet}{\bullet}{\bullet}{\bullet}};
\path (0,-10)   node (3)  %
    {\arbst{\bullet}{\bullet}{\bullet}{\bullet}{\bullet}};
\path (15,-10)   node (4)  %
    {\arbsq{\bullet}{\bullet}{\bullet}{\bullet}{\bullet}};
\path (-15,-20)   node (5) %
    {\arbsc{\bullet}{\bullet}{\bullet}{\bullet}};
\path (0,-20)   node (6) %
    {\arbss{\bullet}{\bullet}{\bullet}{\bullet}};
\path (15,-20)   node (7) %
    {\arbsse{\bullet}{\bullet}{\bullet}{\bullet}};
\path (0,-30)   node (8) %
    {\arbsh{\bullet}{\bullet}{\bullet}};

\path[every node/.style={font=\sffamily\small}] ;
\draw[->] (1) -- (2) ;
\draw[->] (1) -- (3) ;
\draw[->] (1) -- (4) ;
\draw[->] (2) -- (5) ;
\draw[->] (2) -- (6) ;
\draw[->] (3) -- (5) ;
\draw[->] (3) -- (7) ;
\draw[->] (4) -- (6) ;
\draw[->] (4) -- (7) ;
\draw[->] (5) -- (8) ;
\draw[->] (6) -- (8) ;
\draw[->] (7) -- (8) ;
\end{tikzpicture}
%\begin{tikzpicture}[->,>=stealth',shorten >=1pt,auto,node distance=5cm,
%                    thick,main
%node/.style={}]
%
%  \node[main node] (1) %
%    {\arbsu{\bullet}{\bullet}{\bullet}{\bullet}{\bullet}{\bullet}};
%  \node[main node] (2) [below left of=1] %
%    {\arbsd{\bullet}{\bullet}{\bullet}{\bullet}{\bullet}};
%  \node[main node] (3) [below right of=2] %
%    {\arbsq{\bullet}{\bullet}{\bullet}{\bullet}};
%  \node[main node] (4) [below right of=1] %
%    {\arbst{\bullet}{\bullet}{\bullet}{\bullet}{\bullet}};
%
%
%  \path[every node/.style={font=\sffamily\small}]
%    (1) edge node [left] {} (4)
%        edge node [right] {} (2)
%    (2) edge node [right] {} (3)
%    (4) edge node [left] {} (3)
%  ;
%\end{tikzpicture}
\\
which are respectively, reading the diagram by rows, 
$T(201100011000)$, $T(110200011000)$, $T(110110002000)$,
$T(20200011000)$, $T(20110002000)$, $T(11020002000)$,
and $T(2020002000)$.

}
\end{example}
%%%%%%%%%%%%%%%%%%%%%%%%%%%%%%%%%%%%%%%%%%%%%%%%%%%%%%%%%%%%%%%%%%%%%%%%%%%%%%%
%%%%%%%%%%%%%%%%%%%%%%%%%%%%%%%%%%%%%%%%%%%%%%%%%%%%%%%%%%%%%%%%%%%%%%%%%%%%%%%
%%%%%%%%%%%%%%%%%%%%%%%%%%%%%%%%%%%%%%%%%%%%%%%%%%%%%%%%%%%%%%%%%%%%%%%%%%%%%%%
\section{A binary tree expansion}

%%%%%%%%%%%%%%%%%%%%%%%%%%%%%%%%%%%%%%%%%%%%%%%%%%%%%%%%%%%%%%%%%%%%%%%%%%%%%%%
\subsection{A bilinear map on $\Sym$}

The preceding remark (Note \ref{note:bin}) suggests that, as in the case
of Lagrange inversion, the conjugacy equation can be cast as a quadratic
fixed point problem. This is easily done at the level of noncommutative
symmetric functions.

Let $\Omega$ be the linear operator on $\Sym$ introduced in \cite{NT1,NTLag},
and defined by
\begin{equation}\label{def:omeg}
\Omega S^{(i_1,\ldots,i_r)} = S^{(i_1+1,i_2,\ldots,i_r)}\,.
\end{equation}
Writing
\begin{align}
g(qtA)-g(tA)
     &= \sum_{n\ge 1}t^nS_n(A) g(tA)^{n+1}\\
     &= \left(\sum_{n\ge 1}t^n S_n(A)g(tA)^{n}\right)g(tA)\\
     &= tS_1(A)g(tA)^2+t\left(\sum_{n\ge 2}t^{n-1}S_n(A) g(tA)^{n}\right)g(tA)\\
     &= tS_1(A)g(tA)^2+t\Omega\left[\sum_{n\ge 1}t^n S_n(A)g(tA)^{n+1}\right]g(tA)\\
     &= tS_1(A)g(tA)^2 + t\Omega[g(qtA)-g(tA)]g(tA)
\end{align}
we see that $g$ is the unique solution of the quadratic functional equation
\begin{equation}
g(tA)= 1 + \dashint_0^t\left( S_1(A)g(sA)+t\Omega\Delta_qg(sA)\right)g(sA)d_qs
\end{equation}
of the form
\begin{equation}\label{eq:Bq}
g = 1+\B_q(g,g).
\end{equation}

The bilinear map $\B_q$ has a simple expression in the complete basis. For
two compositions $I\vDash i$ and $J\vDash j$,
\begin{align}
\label{BonS}
\B_q(S^I,S^J)
 &= \dashint_0^1\left( S_1S^IS^Js^{i+j}+\Omega S^I(q^i-1)s^iS^Js^j)\right)d_qs \\
 &= \frac{S^{1IJ}+q_i\Omega S^{IJ}}{q_{i+j+1}}\,.
\end{align}
It follows that in the ribbon basis
\begin{equation}\label{BonR}
\B_q(R_I,R_J)
 = \frac{(R_{1I}+q^iR_{1\triangleright I})R_J}{q_{i+j+1}}\,.
\end{equation}
As a consequence, the coefficients of $\tg_n$ on the ribbon basis are in
$\N[q]$.
A combinatorial interpretation will be provided below.

For example,  with $g_0=1$, we have,
\begin{equation}
\begin{split}
g_1 &= \B_q(g_0,g_0) \\
g_2 &= \B_q(g_0,g_1) + \B_q(g_1,g_0) \\
g_3 &= \B_q(g_0,g_2) + \B_q(g_1,g_1) + \B_q(g_2,g_0) \\
\end{split}
\end{equation}
so that one gets
\begin{equation}
\begin{split}
\label{ex-bq}
\tg_1 &= q_1 \B_q(\tg_0,\tg_0) \\
\tg_2 &= q_2 (\B_q(\tg_0,\tg_1) + \B_q(\tg_1,\tg_0)) \\
\tg_3 &= q_3 (\B_q(\tg_0,\tg_2) + \frac{q_2}{q_1}\B_q(\tg_1,\tg_1)
         + \B_q(\tg_2,\tg_0)) \\
\end{split}
\end{equation}
and
\begin{equation}
\begin{split}
\tg_1 &= S_{1} \\
\tg_2 &= (S^{11}) +(S^{11}+q_1S_2) \\
\tg_3 &= (2S^{111}+q_1S^{12}) + \frac{q_2}{q_1}(S^{111}+q_1S^{21}) 
         + (2S^{111}+q_1S^{12}+2q_2S^{21}+q_1q_2S_3) \\
      &= q_1q_2S_{3} + 3q_2S^{21} + 2q_1S^{12} + \left(\frac{q_2}{q_1}+4\right)S^{111},
\end{split}
\end{equation}
which coincides with \eqref{first-tg}.

%%%%%%%%%%%%%%%%%%%%%%%%%%%%%%%%%%%%%%%%%%%%%%%%%%%%%%%%%%%%%%%%%%%%%%%%%%%%%%%

%%%%%%%%%%%%%%%%%%%%%%%%%%%%%%%%%%%%%%%%%%%%%%%%%%%%%%%%%%%%%%%%%%%%%%%%%%%%%%%
\subsection{Triduplicial expansion}

These equations can be lifted to Schr\"oder trees, by setting as above
$g_0=S_0$ and using \eqref{BonS} without modification. We recover then the
same expressions for $g_n$ as in the previous section. 

Indeed, start again with~\eqref{ex-bq} and $\tg_0=g_0=S_0$.
We then have
\begin{equation}\label{ex-tridexp}
\begin{split}
\tg_1 &= S^{100}, \\
\tg_2 &= (S^{10100}) +(S^{11000}+q_1S^{2000}), \\
\tg_3 &= (S^{1010100}+S^{1011000}+q_1S^{102000})
       + \frac{q_2}{q_1}(S^{1100100}+q_1S^{200100})  \\
      &+  (S^{1101000}+S^{1110000}+q_1S^{120000}
       + q_2S^{201000}+q_2S^{210000}+q_1q_2S^{30000}), \\
\end{split}
\end{equation}
which again gives back \eqref{first-tg0}.

\bigskip
Now, \eqref{BonR} can be lifted in another way. Indeed, in the same way as
Lagrange inversion is directly related to Catalan numbers (in the guise of nondecreasing
parking functions) and to  the free duplicial algebra on one
generator $\CQSym$ \cite{NTDup}, 
we find here the little Schr\"oder numbers which are related  to the
free triduplicial algebra on one generator (defined in \cite{NTDup}),
into which $\CQSym$ is naturally embedded.

More precisely, the co-hypoplactic subalgebra of $\PQSym$, denoted by $\SQSym$ in \cite{NT1},
spanned by hypoplactic classes of parking functions (parking quasi-ribbons),  has been
identified in \cite{NTDup} as the free triduplicial algebra on one generator.
Its graded dimension is given by the number of Schr\"oder trees, but no natural
bijection between these trees and parking quasi-ribbons had been known up to now.
However, this algebra has a basis $\PS_\alpha$ which is mapped to the ribbon
basis $R_I$ by a Hopf algebra morphism $\chi$. This suggests that the
operations on ribbons involved in \eqref{BonR} might be the image under $\chi$
of triduplicial operations on parking quasi-ribbons
and that an analogue of the $S$-basis should exist in $\SQSym$. We shall see
that this is indeed the case, by means of bijections between  three families of Schr\"oder
objects: parking quasi-ribbons, Schr\"oder trees, and Schr\"oder pseudocompositions.
These bijections  will allow  to transport the triduplicial structure on the latter objects.

Recall from \cite{NT1} that hypoplactic classes of parking functions are
represented as parking quasi-ribbons, or segmented nondecreasing parking
functions, \emph{i.e.}, nondecreasing parking functions with bars allowed
between different values, for example
\begin{equation}
\label{ex12}
\{1\}, \qquad\qquad \{11,\, 12,\, 1\sep2\},
\end{equation} 
\begin{equation}
\label{ex3}
\{111,\ \ 112,\ \ 11\sep2, \ \ 113,\ \ 11\sep3,\ \ 122,\ \ 1\sep22,\ \
123,\ \ 1\sep23,\ \ 12\sep3,\ \ 1\sep2\sep3 \}.
\end{equation} 

With a parking quasi-ribbon $\alpha$, we associate the elements
\begin{equation}
\PS_{\alpha} := \sum_{\overline\park=\alpha} {\F_\park},
\end{equation}
where $\overline \park$ denotes the hypoplactic class of $\park$.
For example,
\begin{equation}
\PS_{11|3} = \F_{131}+\F_{311}\,, \qquad
\PS_{113} = \F_{113}.
\end{equation}
The product formula in this basis is
\begin{equation}
\PS_{\alpha} \PS_{\beta} = \PS_{\alpha|\beta'} +
\PS_{\alpha\cdot\beta'}
\end{equation}
where $\beta'=\beta[|\alpha|]$ ({\it i.e.}, the word
formed by the entries of $\beta$ shifted by the length of $ \alpha$),
and the dot denotes concatenation.
The triduplicial operations on parking quasi-ribbons are defined
by~\cite{NTDup}
\begin{align}
\alpha \prec\beta &= \alpha\cdot \beta[\max(\alpha)-1],\\
\alpha\circ\beta &= \alpha\sep \beta[|\alpha|],\\
\alpha\succ\beta &= \alpha\cdot \beta[|\alpha|].
\end{align}
One easily checks that they satisfy the seven triduplicial relations
\begin{equation}
\begin{split}
(x\prec y)\prec z  &=  x\prec(y\prec z) \\
(x\circ y)\circ z  &=  x\circ(y\circ z) \\
(x\succ y)\succ z  &=  x\succ(y\succ z) \\
(x\succ y)\prec z  &=  x\succ(y\prec z) \\
(x\circ y)\prec z  &=  x\circ(y\prec z) \\
(x\succ y)\circ z  &=  x\succ(y\circ z) \\
(x\circ y)\succ z  &=  x\circ(y\succ z).
\end{split}
\end{equation}

In order to define the triduplicial operations on Schr\"oder
pseudocompositions, we first need a bijection, which will be described
below.

%%%%%%%%%%%%%%%%%%%%%%%%%%%%%%%%%%%%%%%%%%%%%%%%%%%%%%%%%%%%%%%%%%%%%%%%%%%%%%
\subsection{A bijection between parking quasiribbons and Schr\"oder 
trees}

The bijection between Schr\"oder
pseudocompositions and Schr\"oder trees is trivial, as
it is essentially the Polish notation for the tree.
The difficult point is the correspondence between trees
and parking quasi-ribbons.

Among all Schr\"oder trees, we have binary trees,  and among parking
quasi-ribbons, we have parking quasi-ribbons without bars, that are
nondecreasing parking functions. Both are counted by Catalan numbers.

We shall first describe the bijection from binary trees to parking quasi-ribbons
without bars. Its extension to all Schr\"oder trees will then be straigthforward.
Let $\phi$ be this bijection. It is recursively defined as follows. Set $\phi(\emptyset)=\emptyset$
and $\phi(\bullet)=1$.

Given a tree $T$ with left and right subtrees respectively $T_1$ and
$T_2$, we have
\begin{equation}
\phi(T) = \phi(T_2) \cdot (\phi(T_1)[\max(\phi(T_2))-1])\cdot (|T_1|+\max(\phi(T_2)),
\end{equation}
with the convention that, if $T_2$ is empty, $\max(\phi(T_2))=1$, and the dot denotes concatenation.
This operation can also be described as collecting the vertices of $T$
recursively by visiting first its right subtree, then its left subtree and finally  its root,
 with the rules
that a leaf takes the value of the last visited vertex ($1$ if there were none)
and an internal vertex gets as value the size of its left subtree added to the
value of its right son (added to $1$ if there is no right son).

\begin{example}{\rm
We have
\begin{equation}
\label{ex-bijab}
{\xymatrix@C=2mm@R=2mm{
        *{} & *{} & {8}\ar@{-}[drr]\ar@{-}[dll] \\
        {6}\ar@{-}[dr]  & *{} & *{} & *{} & {5}\ar@{-}[dr]\ar@{-}[dl] \\
        *{}  & {6}\ar@{-}[dl] & *{} & {4} & *{} & 4\ar@{-}[dl] \\
        {5}  & *{} & *{} & *{} & {2}\ar@{-}[dr]\ar@{-}[dl] & *{} & {} \\
        *{}  & *{} & *{} & {1} & *{} & {1} & {} 
      }}
\quad \longrightarrow 1124455668.
\end{equation}
}
\end{example}

Let us now extend $\phi$ to all Schr\"oder trees.
First, Schr\"oder trees are in bijection with binary trees with two-colored
left edges: if an internal node $s$ has more than two children
with corresponding subtrees $T_1,T_2,\dots,T_r$, draw $r-1$ left edges (of
the second color) from $s$ and attach to the new $r$ nodes the $r$ subtrees in
order, as in a binary tree:
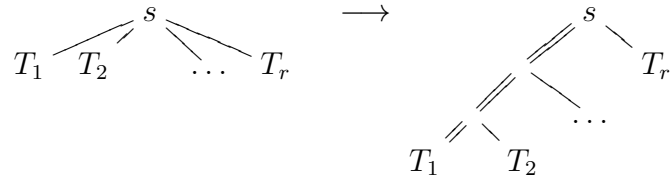
\begin{equation*}
{\xymatrix@C=2mm@R=2mm{
        *{} & *{} & {s}\ar@{-}[drr]\ar@{-}[dl]\ar@{-}[dr]\ar@{-}[dll] \\
        {T_1} & {T_2} & *{} & {\dots} & {T_r}
      }}
\quad
\longrightarrow
{\xymatrix@C=2mm@R=2mm{
        *{} & *{} & *{} & {s}\ar@{=}[dl]\ar@{-}[dr] \\
        *{} & *{} & {}\ar@{=}[dl]\ar@{-}[dr] & *{} & T_r \\
        *{} & {}\ar@{=}[dl]\ar@{-}[dr] & *{} & \dots \\
        {T_1} & *{} & {T_2}
      }}
\end{equation*}
This amounts to reverting the contraction process described in Note \ref{note:bin}.

Having computed this tree, send it with the previous bijection to a nondecreasing
parking function, and insert a bar between two letters if they are separated by a
left branch of the second color.

\begin{example}{\rm
The continuation of \eqref{ex-bijab} is
\begin{equation}
{\xymatrix@C=2mm@R=2mm{
        *{} & *{} & {8}\ar@{-}[drr]\ar@{-}[dll] \\
        {6}\ar@{-}[dr]  & *{} & *{} & *{} & {5}\ar@{-}[dr]\ar@{-}[dl] \\
        *{}  & {6}\ar@{-}[dl] & *{} & {4} & *{} & 4\ar@{-}[dl] \\
        {5}  & *{} & *{} & *{} & {2}\ar@{-}[dr]\ar@{-}[dl] & *{} & {} \\
        *{}  & *{} & *{} & {1} & *{} & {1} & {} 
      }}
\quad\longrightarrow
{\xymatrix@C=2mm@R=2mm{
        *{} & *{} & {8}\ar@{=}[dll]\ar@{-}[drr] \\
        {6}\ar@{-}[dr]  & *{} & *{} & *{} & {5}\ar@{=}[dl]\ar@{-}[dr] \\
        *{}  & {6}\ar@{-}[dl] & *{} & {4} & *{} & 4\ar@{-}[dl] \\
        {5}  & *{} & *{} & *{} & {2}\ar@{=}[dl]\ar@{-}[dr] & *{} & {} \\
        *{}  & *{} & *{} & {1} & *{} & {1} & {} 
      }}
\quad \longrightarrow 11|244|5566|8.
\end{equation}
}
\end{example}
\begin{theorem}
The previous algorithm provides a bijection between Schr\"oder trees and parking
quasi-ribbons.
\end{theorem}

Before proving the theorem, let us describe the inverse bijection  $\psi=\phi^{-1}$
from nondecreasing parking functions to binary trees.
Set $\psi(1)=\emptyset$.
Let $p=a_1\dots a_r$ be a nondecreasing parking function. Let
$w=w_1\dots w_{r-1}$ be the word such that $w_k=a_k+r-1-k$. 
Let $\ell$ be greatest index
of $w$ such that $a_{\ell}=a_{\ell+1}$ and $w_\ell=a_r$ (as $a_0=1$, if $\ell$ does not
exist,  set $\ell=0$). Then compute recursively the images of
$a_1\dots a_\ell$ as the right subtree of the root, and of
$(a_{\ell+1}\dots a_{r-1})[a_\ell-1]$ as the left subtree of the root.

\begin{example}{\rm
Consider all nondecreasing parking functions $112445566X$, where $X$ is bound
by the constraint of being a nondecreasing parking function, so that $X\in \{6,7,8,9,10\}$.
The word $w$ is the same for all those parking functions, namely
${\bf9}88{\bf9}8{\bf8}7{\bf7}6$, where we write in boldface the 
$w_k$ such that
$a_k=a_{k+1}$. Note that the $6$ can be bold or not
depending on the value of $X$: if $X$ is $6$, it is indeed in bold.
Now, the index $\ell$ is well-defined in all the examples, so that we can
separate the word and apply it recursively:
\begin{equation}
\label{ex-bijabb}
\begin{split}
1124455666    & \longrightarrow_{\ell=9} (\emptyset, 112445566) \\
1124455667    & \longrightarrow_{\ell=8} (1, 11244556) \\
1124455668    & \longrightarrow_{\ell=6} (122, 112445) \\
1124455669    & \longrightarrow_{\ell=4} (12233, 1124) \\
112445566\,10 & \longrightarrow_{\ell=0} (112445566, \emptyset) \\
\end{split}
\end{equation}
}
\end{example}
\Proof
Let us now prove that $\phi$ is indeed a bijection and that its inverse is
$\psi$ as claimed.

First, the values of $\phi$  are clearly nondecreasing parking functions.

For $\psi$, the crucial point  is to prove that it is
well-defined (see \eqref{ex-bijabb} for an illustration).

Given a nondecreasing parking function $p=a_1\dots a_r$, the allowed values for $a_r$ 
are in the interval $[a_{r-1},r]$. 
And this interval corresponds precisely
to the values taken by the subword of $w$ obtained by selecting  the indices $i$
such that $a_i=a_{i+1}$. 
Indeed, any of these values belong to this interval,
since all values of $w$ do. 
Conversely, a direct induction on the 
length of $p$ implies the result, since the only question concerns the index
$r-1$ which is considered in the subword of $w$ iff $a_{r-1}=a_r$. 
Finally,
if one splits $p$ after the rightmost occurrence $w_\ell$ of such a value,
both the prefix of $p$ and its suffix $a_{\ell+1}\dots a_{r-1}$ are
parking functions: it is obvious for the prefix and is easy for the suffix,
since we considered the \emph{rightmost occurrence}.
This occurrence
has only strictly smaller values to its left ($w_i-w_{i+1}$ can be at most
$1$), so that the shifted suffix is a parking function.
Now, by definition, the values of $\psi$ are binary trees, so that at this point, we
have  maps going from each set to the other. Let us now see
why they are inverses of each other.

Both maps are recursive, so we just need to prove that they are inverse of each
other on the first step. Let $p=a_1\dots a_r$ be a nondecreasing parking function which is the
image under $\phi$ of a binary tree $T$ having $T_1$ and $T_2$ as
left and right subtrees. Since $a_r$ is  the sum of the
size of $T_1$ and of the maximal value of $\phi(T_2)$, $a_r$ corresponds to the
value $w_\ell$ in the word $w$ where $\ell$ is the size of $T_2$ and
$r-\ell-1$ is the size of $T_1$. Now, this value $\ell$ necessarily satisfies
$a_\ell=a_{\ell+1}$, since $a_\ell$ is the maximal value $M$
of $T_2$ and $a_{\ell+1}=1+M-1$. Finally, among all indices $k$ satisfying
$a_k=a_{k+1}$ and $w_k=w_\ell$, $\ell$ is the only one such that the suffix
$a_{\ell+1}\dots a_{r-1}$ is a parking function, since any other occurrence in $w$
has one equal value to its right, which contradicts the fact of being a
nondecreasing parking function.

Since $\phi$ and $\psi$ both have the right image sets and
$\psi\circ\phi$ is the identity map on binary trees, they both
are bijections, inverses of each other.

\medskip
Finally, let us see why the extension of $\phi$ to Schr\"oder trees is a
bijection. First, all left branches relate numbers that cannot be equal, so that
separations on nondecreasing parking functions are made between non equal
letters, which is the required condition about parking quasi-ribbons.
The converse is also true: the number of left branches in a binary tree
$T$ is equal to the number of different letters plus one in $\phi(T)$.
So the map from binary trees with two-colored left branches to parking
quasi-ribbons is a bijection and the composition of both bijections through
the middle object of binary trees with two-colored left branches is still a
bijection. 
\qed

%%%%%%%%%%%%%%%%%%%%%%%%%%%%%%%%%%%%%%%%%%%%%%%%%%%%%%%%%%%%%%%%%%%%%%%%%%%%%%
\subsection{Triduplicial operations on Schr\"oder pseudocompositions}

Now that we have a bijection between parking quasiribbons and Schr\"oder
pseudocompositions, we can translate the triduplicial operations 
initially defined on parking quasiribbons to Schr\"oder pseudocompositions.

\begin{definition}
Let $I$ and $J$ be two Schr\"oder pseudocompositions.
Define $J'$ such that $J=J'0^{m}$ and $J'$ does not end by a 0.

Then
\begin{align}
I\prec J &= J\triangleright I = J'0^{m-1}.I\\
I\circ J &= J'\triangleright I\cdot 0^{m-1}\\
I\succ J &= J'\cdot I\cdot 0^{m-1}\,.
\end{align}
\end{definition}

\begin{example}{\rm
Denoting by $a$ the parking quasi-ribbon $1$ and by $x$ the
pseudocomposition $100$,
\begin{align}
a\prec a = 11\quad   & x\prec x = 10100\\
a\circ a = 1|2\quad  & x\circ x = 2000\\
a\succ a = 12\quad   & x\succ x = 11000,
\end{align}
which coincide with the bijection described
in the previous section.
}
\end{example}

\begin{theorem}
The operations above endow the set of Schr\"oder pseudocompositions with the
structure of a triduplicial algebra, freely generated by $x=100$.
% This induces
%a bijection between parking quasi-ribbons and  Schr\"oder pseudocompositions.
\end{theorem}

\Proof
This is a direct consequence of the translation of the
triduplicial operations on Schr\"oder pseudocompositions.

We shall prove it for each rule.
Operation $\prec$ on parking quasi-ribbons is defined by
$\alpha \prec\beta = \alpha\cdot \beta[\max(\alpha)-1]$ and via the bijection
 $\psi$ extended to parking quasi-ribbons,  it
corresponds to glueing the image of $\alpha$ to the rightmost leaf of $\beta$
so that, on Schr\"oder pseudocompositions, one obtains $J'.0^{m-1}.I$.

Operations $\circ$ and $\succ$ on parking quasi-ribbons are defined by
$\alpha\circ\beta = \alpha\sep \beta[|\alpha|]$
and
$\alpha\succ\beta = \alpha\cdot \beta[|\alpha|]$
and via the bijection $\psi$ extended to parking quasi-ribbons, 
it corresponds to putting the image of $\alpha$ as the left child of
the rightmost internal node labelled $1$ of the image of $\beta$ (which is
also the last visited internal node of the tree in Polish notation) with an edge
of the natural color (for $\succ$) or of the second color (for $\circ$).
The translation on Schr\"oder pseudocompositions is straightforward.
\qed

Define now an order $\le$ on parking quasiribbons by the cover relation
\begin{equation}
\beta \gtrdot \alpha\quad\text{if $\alpha=uv$, $\beta=u|v'$}
\end{equation}
with $v'=v$ if the last letter of $u$ is smaller than the first letter of $v$,
and $v'=v[1]$ otherwise.

For example, the predecessors of $11|23$ are $11|2|3$ and $1|2|34$.

With this order, we can define a basis $\S^\alpha$ by
\begin{equation}
\S^\alpha=\sum_{\alpha\le\beta}\PS_\beta\,.
\end{equation}

For example,
\begin{equation}
\S^1=\PS_1,\ \S^{11}= \PS_{11}+\PS_{1|2},\
\S^{12}=\PS_{12}+\PS_{1|2},\ \S^{1|2}=\PS_{1|2}
\end{equation}
and
\begin{equation}
\S^{11|23}=\PS_{11|23}+\PS_{1|2|34}+\PS_{11|2|3}+\PS_{1|2|3|4}.
\end{equation}

The Hopf epimorphism $\chi:\ \SQSym\rightarrow\Sym$ is defined by
\begin{equation}
\chi(\PS_\alpha) = R_{I^\sim}
\end{equation}
where $I$ is the bar composition of $\alpha$ whose parts are the lengths of
the factors between the bars, {\it e.g.}, for $\alpha=111|24|5$, $I=321$.

It is then clear that
\begin{equation}
\chi(\S^\alpha)=S^{I^\sim}.
\end{equation}

For $U_m\in\SQSym_m$ and $V_n\in\SQSym_n$, define
\begin{equation}
B_q(U_m,V_n) = q_{m+n+1}\B_q(U_m,V_n),
\end{equation}
where $\B_q$ is defined after Eq. \eqref{eq:Bq}.
Then, in the bases $\PS$ and $\S$ with both indexations:

\begin{align}
B_q(\S^\alpha,\S^\beta) &=
 q_{|\alpha|}\S^{\beta\prec(\alpha\circ 1)}+ \S^{\beta\prec(\alpha\succ 1)} \\
B_q(\S^I,\S^J) &=
 q_{|I|}\Omega\S^{IJ}+\S^{1IJ} \\
B_q(\PS_\alpha,\PS_\beta) &=
 q^{|\alpha|}\PS_{\beta\prec(\alpha\circ 1)}+\PS_{\beta\prec(\alpha\succ 1)}
 +q^{|\alpha|}\PS_{\beta\circ\alpha\circ 1}+\PS_{\beta\circ\alpha\succ 1}\\
B_q(\PS_I,\PS_J) &=
 q^{|I|}\PS_{\Omega IJ}+\PS_{1IJ}
 +q^{|I|}P_{1\triangleright J'\triangleright I0^m}
 +\PS_{1J'\triangleright I0^m}
\end{align}
where, as above, $J=J'0^m$ and $J'$ does not end by a $0$.

%%%%%%%%%%%%%%%%%%%%%%%%%%%%%%%%%%%%%%%%%%%%%%%%%%%%%%%%%%%%%%%%%%%%%%%%%%%%%%%%

For example, one can recover the computation of $\tg_3$ in
\eqref{ex-tridexp}: start from the expression of $\tg_1$ and $\tg_2$ in this
same equation and then compute $\tg_3$ according to Eq.~\eqref{ex-bq}:

The first term is $B_q(\tg_0,\tg_2)$:
\begin{equation}
\begin{split}
B_q(S_0,S^{10100}) &= q_0 S^{110100} + S^{1010100} = S^{1010100}, \\
B_q(S_0,S^{11000}) &= q_0 S^{111000} + S^{1011000} = S^{1011000}, \\
B_q(S_0,S^{2000} ) &= q_0 S^{12000} + S^{102000} = S^{102000}. \\
\end{split}
\end{equation}
Note that some $S^I$ here are not indexed by Schr\"oder pseudocompositions, but
these terms eventually disappear
as their coefficient is $q_0=q^0-1=0$.
The second term is $B_q(\tg_1,\tg_1)$:
\begin{equation}
B_q(S^{100}, S^{100}) = q_1 S^{200100} + S^{1100100}. \\
\end{equation}
The third term is $\B_q(\tg_2,\tg_0)$:
\begin{equation}
\begin{split}
B_q(S^{10100}, S_0) &= q_2 S^{201000} + S^{1101000}, \\
B_q(S^{11000}, S_0) &= q_2 S^{210000} + S^{1110000}, \\
B_q(S^{2000},  S_0) &= q_2 S^{30000}  + S^{120000}, 
\end{split}
\end{equation}
so that we recover \eqref{ex-tridexp}.
%%%%%%%%%%%%%%%%%%%%%%%%%%%%%%%%%%%%%%%%%%%%%%%%%%%%%%%%%%%%%%%%%%%%%%%%%%%%%%%%

%%%%%%%%%%%%%%%%%%%%%%%%%%%%%%%%%%%%%%%%%%%%%%%%%%%%%%%%%%%%%%%%%%%%%%%%%%%%%%%%
\section{Expansion on the ribbon basis}\label{sec:rib}

The expression of $g$ in $\Sym$ is recovered by setting $S_0=1$. 
As in the case of the Lagrange series, it is interesting to expand $g$ 
on the ribbon basis. As we have already seen before, the first terms are
\begin{align}
\tg_1 &= R_1,\\
\tg_2 &= (1+q)R_2 + 2R_{11},\\
\tg_3 &= (1+q)(1+q+q^2)R_3+ (2+q+3q^2)R_{21}+3(1+q)R_{12}+(5+q)R_{111}. 
\end{align} 

We can observe that each coefficient is a $q$-analogue of $n!$.
We shall now prove this fact, and describe the relevant statistics on permutations.

For a pseudo-composition $I$, let $\hat I$ be the ordinary composition
obtained by removing the zero entries.

For a binary tree $t =T(I)$, set
\begin{equation}
P_t = \sum_{J\preceq I}m_J(q)S^J
\end{equation}
For a pseudo-composition $J$ encoding a tree $T(J)$, let
\begin{equation}
d_J = \sum_{v\in T(J)}(\phi(v)-1)
\end{equation}
Then, the coefficient of $R_{K}$ in $(q)_nP_t$ is equal to
\begin{equation}
(q)_nm_I(q)q^{d_I-d_J}
\end{equation}
where $J$ is the coarsest anti-refinement of $I$ such that $K\le \hat J$.
Indeed, $R_K$ will then occur in all the refinements of $J$, 
and if $I'$ is such a refinement, then,
\begin{equation}
m_{I'}(q)=m_I(q)\prod_{v\in C(I,I')} q_{\phi(v)-1}
\end{equation} 
where $C(I,I')$ is the set of vertices of $T(I)$ which have been contracted
in $T(I')$.
Thus, factoring
the coefficient $m_I(q)$, we see that $R_K$ picks up a factor 
\begin{equation}
\prod_{i}(q_i+1)^{n_i}=q^{d_I-d_J}\quad \text{if}\quad \frac{m_J(q)}{m_I(q)}=\prod_i q_i^{n_i} 
\end{equation} 
when summing over the Boolean lattice of refinements of $J$.

\begin{example}{\rm
For 
$I= (1101100011000)$ and $K=(51)$, we have $J=(20200011000)$,
and on the picture (extracted from \ref{ex:contracts})\\
\begin{tikzpicture}[->,>=stealth',shorten >=1pt,auto,node distance=5cm,
                    thick,main node/.style={}]

  \node[main node] (1) {\arbsu{6}{3}{2}{2}{1}{1}};
  \node[main node] (2) [below left of=1] {\arbsd{6}{2}{2}{1}{1}};
  \node[main node] (3) [below right of=1] {\arbst{6}{3}{2}{2}{1}};
  \node[main node] (4) [below right of=2] {\arbsc{6}{2}{2}{1}};

  \path[every node/.style={font=\sffamily\small}]
    (1) edge node [left] {} (3)
        edge node [right] {} (2)
    (2) edge node [right] {} (4)
    (3) edge node [left] {} (4)
  ;
\end{tikzpicture}

%%%%%%%%%%%
%\begin{tikzpicture}[->,>=stealth',shorten >=1pt,auto,node distance=5cm,
%                    thick,main
%node/.style={}]
%
%  \node[main node] (1) {\arbsu{6}{3}{2}{2}{1}{1}};
%  \node[main node] (2) [below left of=1] {\arbsd{6}{2}{2}{1}{1}};
%  \node[main node] (3) [below right of=2] {\arbsq{6}{2}{2}{1}};
%  \node[main node] (4) [below right of=1] {\arbst{6}{3}{2}{2}{1}};
%
%  \path[every node/.style={font=\sffamily\small}]
%    (1) edge node [left] {} (4)
%        edge node [right] {} (2)
%    (2) edge node [right] {} (3)
%    (4) edge node [left] {} (3)
%  ;
%\end{tikzpicture}
we can read that the coefficient of $R_{51}$ in the projection of
$(q)_6P_{T(1101100011000)}$ on $\Sym$ is
\begin{equation}
\frac{(q)_6}{q_6q_3q_2^2q_1 ^2}q^4 = \frac{q_5q_4}{q_2q_1}q^4 = {q}^{10}+2{q}^{8}+{q}^{9}+2{q}^{7}+2{q}^{6}+{q}^{5}+{q}^{4}.
\end{equation}
}
\end{example}
For a pseudo-composition $I$, let $I^\sharp$ be its coarsest anti-refinement. The polynomial
\begin{equation}\label{eq:qh}
(q)_nm_I(q)q^{d_I-d_{I^\sharp}}
\end{equation}
is, up to left-right symmetry, the $q$-hook-length formula for the binary tree $T(I)$ \cite{BW,HNT} (with
its leaves removed).
Indeed, for a vertex $v$ of a binary tree $t=T(I)$, $d_I-d_{I^\sharp}$ coincides with the number of
internal nodes of its left subtree.

\begin{example}{\rm
Continuing with $I=(110110011000)$, $I^\sharp=(2020002000)$,
\begin{equation}
T(I)=\arbsu{6}{3}{2}{2}{1}{1}
\qquad
T(I^\sharp)=\arbsh{6}{2}{2}
\end{equation}
and \eqref{eq:qh} yields
\begin{equation}
\frac{(q)_6}{q_6q_3q_2^2q_1^2}q^{6+3+2+1+2+1-(6+2+2)}=\frac{q_5q_4}{q_2q_1}q^{3+1+1} =\frac{q_5q_4}{q_2q_1}q^{5}
\end{equation}
which is the $q$-hook-length formula for the left-right flip
\def\arbeq#1#2#3#4#5#6{
%\newdimen\vcadre\vcadre=0.01cm % marges verticales de la boite
%\newdimen\hcadre\hcadre=0.01cm % marges horizontales de la boite
\xymatrix@R=0.1cm@C=2mm{
& & {\GrTeXBox{#1}}\arx1[lld]\arx1[rrd]\\
% ligne 2
{\GrTeXBox{#2}}\arx1[dr] &&&&
  {\GrTeXBox{#3}}\arx1[dl] \\
% ligne 3
& {\GrTeXBox{#4}} && {\GrTeXBox{#5}}\arx1[dr] \\
% ligne 4
&&&& {\GrTeXBox{#6}} \\
}
} 
\begin{equation}
\arbeq{6}{2}{3}{1}{2}{1}
\qquad
\arbeq{3}{1}{0}{0}{1}{0}
\end{equation}
of the skeleton of $T(I)$. Indeed, the power of $q$ is $3+1+1=5$ (given by the sum of the sizes of the right subtrees), and the denominator
is $q_6q_2q_1q_3q_2q_1$ (recording the sizes of all the subtrees).
}
\end{example}
Thus, up to a fixed power of $q$, this expression enumerates by number of non-inversions the permutations
in the sylvester class labelled by the binary tree $T(I)$. 

Summing over all binary trees, we see that each coefficient $c_K(q)$
of $R_K$ in $\tg_n$ is indeed a $q$-analogue of $n!$.

Translating these results at the level of permutations yields the
following description of the expansion of $\tg_n$ on the ribbon basis of
$\NCSF$:

\begin{definition}
Let $\sigma$ be a permutation and let $\alpha$ be the top of its sylvester
class, that is, the permutation with the smallest number of inversions in its
sylvester class. Let $I=(i_1,\dots,i_r)$ be a composition and let $D$ be the
descent  set of the conjugate $\tilde I$.
\def\invs{{\rm inv}}
Define $C_I(\sigma)$ as
\begin{equation}
q^{\invs(\sigma) - \invs(\alpha)}  q^{\invs(\alpha,D)},
\end{equation}
where $\invs(\alpha,D)$ is the number of pairs $(i<j)$ such that
$\alpha_i>\alpha_j$ and $j\in D$.
\end{definition}

\begin{theorem}\label{th:rib}
The coefficient of $R_I$ in the expansion of $\tg_n$ is equal to
\begin{equation}
\sum_{\sigma\in\SG_n} C_I(\sigma).
\end{equation}
\qed
\end{theorem}

For example, here are the tables for $n=3$ and $n=4$ of all coefficients
$C_I$, where permutations are grouped by sylvester classes.
%permutations of $\SG_3$ with the exponents
% values for
%all compositions $I$ of size $3$.

\begin{equation}
\begin{array}{|c|c|c|c|c|}
\hline
\sigma / I & 3 & 21 & 12 & 111 \\
\hline
\hline
123 & 1 & 1 & 1 & 1 \\[3pt]
\hline
132 & q & q & 1 & 1 \\[3pt]
312 & q^2 & q^2 & q & q \\[3pt]
\hline
213 & q & 1 & q & 1 \\[3pt]
\hline
231 & q^2 & q^2 & 1 & 1 \\[3pt]
\hline
321 & q^3 & q^2 & q & 1 \\[3pt]
\hline
\end{array}
\end{equation}

\begin{equation}
\begin{array}{|c|c|c|c|c|c|c|c|c|}
\hline
\sigma / I & 4 & 31 & 22 & 211 & 13 & 121 & 112 & 1111 \\
\hline
\hline
1234 & 1 & 1 & 1 & 1 & 1 & 1 & 1 & 1 \\[3pt]
\hline
1243 & q & q & q & q & 1 & 1 & 1 & 1 \\[3pt]
1423 & q^2 & q^2 & q^2 & q^2 & q & q & q & q \\[3pt]
4123 & q^3 & q^3 & q^3 & q^3 & q^2 & q^2 & q^2 & q^2 \\[3pt]
\hline
1324 & q & q & 1 & 1 & q & q & 1 & 1 \\[3pt]
3124 & q^2 & q^2 & q & q & q^2 & q^2 & q & q \\[3pt]
\hline
1342 & q^2 & q^2 & q^2 & q^2 & 1 & 1 & 1 & 1 \\[3pt]
3142 & q^3 & q^3 & q^3 & q^3 & q & q & q & q \\[3pt]
3412 & q^4 & q^4 & q^4 & q^4 & q^2 & q^2 & q^2 & q^2 \\[3pt]
\hline
1432 & q^3 & q^3 & q^2 & q^2 & q & q & 1 & 1 \\[3pt]
4132 & q^4 & q^4 & q^3 & q^3 & q^2 & q^2 & q & q \\[3pt]
4312 & q^5 & q^5 & q^4 & q^4 & q^3 & q^3 & q^2 & q^2 \\[3pt]
\hline
2134 & q & 1 & q & 1 & q & 1 & q & 1 \\[3pt]
\hline
2143 & q^2 & q & q^2 & q & q & 1 & q & 1 \\[3pt]
2413 & q^3 & q^2 & q^3 & q^2 & q^2 & q & q^2 & q \\[3pt]
4213 & q^4 & q^3 & q^4 & q^3 & q^3 & q^2 & q^3 & q^2 \\[3pt]
\hline
2314 & q^2 & q^2 & 1 & 1 & q^2 & q^2 & 1 & 1 \\[3pt]
\hline
2341 & q^3 & q^3 & q^3 & q^3 & 1 & 1 & 1 & 1 \\[3pt]
\hline
2431 & q^4 & q^4 & q^3 & q^3 & q & q & 1 & 1 \\[3pt]
4231 & q^5 & q^5 & q^4 & q^4 & q^2 & q^2 & q & q \\[3pt]
\hline
3214 & q^3 & q^2 & q & 1 & q^3 & q^2 & q & 1 \\[3pt]
\hline
3241 & q^4 & q^3 & q^4 & q^3 & q & 1 & q & 1 \\[3pt]
\hline
3421 & q^5 & q^5 & q^3 & q^3 & q^2 & q^2 & 1 & 1 \\[3pt]
\hline
4321 & q^6 & q^5 & q^4 & q^3 & q^3 & q^2 & q & 1 \\[3pt]
\hline
\end{array}
\end{equation}

%%%%%%%%%%%%%%%%%%%%%%%%%%%%%%%%%%%%%%%%%%%%%%%%%%%%%%%%%%%%%%%%%%%%%%%%%%%%%%%%
%%%%%%%%%%%%%%%%%%%%%%%%%%%%%%%%%%%%%%%%%%%%%%%%%%%%%%%%%%%%%%%%%%%%%%%%%%%%%%%%
%%%%%%%%%%%%%%%%%%%%%%%%%%%%%%%%%%%%%%%%%%%%%%%%%%%%%%%%%%%%%%%%%%%%%%%%%%%%%%%%
\section{Schr\"oder's equation for the inverse of $h$}\label{schrodeq}

Let now $f=h^{-1}$ where 
$h^{-1}\circ\phi\circ h(z) = qz$
or equivalently 
\begin{equation}\label{eq:schro}
f\circ\phi(w)=qf(w)\ (w=h(z)).
\end{equation}
This is Schr\"oder's equation.
In the noncommutative setting, with again $\phi(z)=qz\sigma_z(A)$, 
it becomes
\begin{equation}
\sum_{k\ge 0}f_kq^{k+1}w^{k+1}\sigma_w((k+1)A) = q\sum_{n\ge 0}f_nw^{n+1}
\end{equation}
which translates into the recurrence relation
\begin{equation}\label{eq:recfn}
f_n = \sum_{k+l=n} q^k f_k S_l((k+1)A).
\end{equation}

\begin{theorem}
Let $L$ be the linear endomorphism of $\Sym$ defined by
\begin{equation}
L(S^I) = S_{i_1}(A)S_{i_2}((i_1+1)A)
         S_{i_3}((i_1+i_2+1)A)\cdots S_{i_r}((i_1+\cdots +i_{r-1}+1)A).
\end{equation}
Then,
\begin{equation}
f_n = L\left( S_n\left(\frac{A}{1-q} \right)\right) =
      \sum_{I\vDash n}
        \frac{q^{\maj(I)}}%
             {(1-q^{i_1})(1-q^{i_1+i_2})\cdots(1-q^{i_1+\cdots+i_r})} 
      L(S^I),
\end{equation}
where $A/(1-q)$ and $\maj$ is defined as in \cite[6.1]{NCSF2}\footnote{$\maj(I)$ 
is the sum of the descents of $I$, {\it i.e.},$\maj(I)=(r-1)i_1+(r-2)i_2+\cdots+i_{r-1}$ if $I=(i_1,\ldots,i_r)$.}.
\end{theorem}

For example,
\begin{equation}
f_1 = \frac1{1-q}S_1(A)
\end{equation}
\begin{equation}
f_2 = \frac{1}{1-q^2}S^2(A)+ \frac{q}{(1-q)(1-q^2)}S_1(A)S_1(2A)
\end{equation}
\begin{equation}
\begin{split}
f_3 =& \frac{1}{1-q^3}S_3(A)+\frac{q^2}{(1-q^2)(1-q^3)}S_2(A)S_1(3A)\\
     &+\frac{q}{(1-q)(1-q^3)}S_1(A)S_2(2A)
      +\frac{q^3}{(1-q)(1-q^2)(1-q^3)}S_1(A)S_1(2A)S_1(3A)
\end{split}
\end{equation}

\Proof
Replacing $S_l((k+1)A)$ by $S_l(A)$ in \eqref{eq:recfn}, we obtain a
recurrence relation satisfied by the expansion on the basis $S^I$  of the
$S_n(A/(1-q))$.
To recover $f_n$ from this expression, we just have to replace each factor
$S_{i_k}$ of $S^I$ by $S_{i_k}((i_1+\cdots +i_{k-1}+1)A)$.
\qed

\begin{note}{\rm The noncommutative symmetric function
\begin{equation}
K_n(q) = (q)_n S_n\left(\frac{A}{1-q} \right) = \sum_{I\vDash n}q^{\maj(I)} R_I
\end{equation}
is the $q$-Klyachko function. For $q$ a primitive $n$th root of 1, it is
mapped to Klyachko's Lie idempotent in the descent algebra of the symmetric group \cite{NCSF1,NCSF2}.
It is also a noncommutative analogue of the Hall-Littlewood function $Q'_{1^n}$ \cite{Hiv}.
}
\end{note}

Naturally, we may also choose to define $f$ by
\begin{equation}\label{eq:recfnb}
f_n = \sum_{k+l=n}q^k S_l((k+1)A)f_k.
\end{equation}
With this choice
\begin{equation}
f = \overline{\left( L\left( \overline
    {\sigma_1\left( \frac{A}{1-q} \right)}\right)\right)}
\end{equation}
and it is what we obtain by mould calculus in the next section (the bar involution
is defined by $\overline{S^I}=S^{\overline I}$, where
$\overline I$ is the mirror composition).

%%%%%%%%%%%%%%%%%%%%%%%%%%%%%%%%%%%%%%%%%%%%%%%%%%%%%%%%%%%%%%%%%%%%%%%%%%%%%%%
%%%%%%%%%%%%%%%%%%%%%%%%%%%%%%%%%%%%%%%%%%%%%%%%%%%%%%%%%%%%%%%%%%%%%%%%%%%%%%%
%%%%%%%%%%%%%%%%%%%%%%%%%%%%%%%%%%%%%%%%%%%%%%%%%%%%%%%%%%%%%%%%%%%%%%%%%%%%%%%
\section{A noncommutative mould expansion}\label{ncmould}
\def\mm{{\mathcal M}}

Ecalle's approach to the linearization equation \eqref{eq:conjug} is to find
a closed expression of the substitution automorphism 
\begin{equation}
H:\ \psi \longmapsto \psi\circ h
\end{equation}
as a differential operator. The idea is to look for an expansion of the form
\begin{equation}
H = \sum_I \mm_I U^I
\end{equation}
where $I=(i_1,\ldots,i_r)$ runs over all compositions,
$U^I=U_{i_1}\cdots U_{i_r}$ as usual, and the $U_n$ are the homogeneous
component of the differential operator 
\begin{equation}
U:\ \psi \longmapsto \psi\circ u,\quad \text{where $\phi(z)=qu(z)$},
\end{equation}
given by the Taylor expansion at $z$ of $\psi(z+ (u(z)-z))$, 
\begin{equation}
U_n = \sum_{I\vDash n} \frac{u^I}{\ell(I)}z^{n+\ell(I)}\partial_z^{\ell(I)}.
\end{equation}
In this setting, the functional equation  \eqref{eq:conjug} reads
\begin{equation}
\label{eq:conjop}
HUM_q =M_qH,\quad\text{where $(M_q\psi)(z)=\psi(qz)$}.
\end{equation}
We have already identified our generic power series $\phi(z)$ with
$qz\sigma_z(X)$, so that $u_n=h_n(X)$. A natural noncommutative analogue
is to set $u_n=S_n(A)$. There is another way to introduce noncommutative
symmetric functions in this problem. The substitution maps $U$ and
$H$, being automorphisms, are grouplike elements in the Hopf algebra
of differential operators. So, it is natural to introduce a second alphabet
$B$ (commuting with $A$) and to identify $U_n$ with $S_n(B)$. The problem
amounts to looking for $H$ as an element of (the completion of)
$\Sym(A)\otimes\Sym(B)$. Let us write it as $H(B)$, regarded as a symmetric
function of $B$ with coefficients in $\Sym(A)$. Then \eqref{eq:conjop}
reads now
\begin{equation}
H(B)\sigma_1(B) = H(qB).
\end{equation}
This is solved by
\begin{equation}
H(B) = \prod_{i\ge 0}^\leftarrow \lambda_{-q^i}(B)
     = \cdots\lambda_{-q^2}(B)\lambda_{-q}(B)\lambda_{-1}(B)
\end{equation}
which may be denoted by
\begin{equation}
H(B) = \sigma_1\left(\frac{B}{q-1}\right)
\end{equation}
Setting $F=H^{-1}$, we have that
\begin{equation}
F(B) = \sigma_1(B)\sigma_q(B)\sigma_{q^2}(B)\cdots 
\end{equation}
is the image of $\sigma_1(A/(1-q))$ (in the sense of \cite{NCSF2}) by the
bar involution $S^I\mapsto S^{\bar I}$, so that

\begin{equation}
F_n = \sum_{I\vDash n}
        \frac{q^{\maj(I)}}{(1-q^{i_1})(1-q^{i_1+i_2})\cdots (1-q^n)}S^{\bar I},
\end{equation}
The function $f(z)$ is obtained by acting on the identity: $f(z)=Fz$. 
This is obtained from
\begin{equation}
U^I z = \overline{L(S^{\bar I}(A))}z^{|I|+1}\,.
\end{equation}
Indeed,
\begin{equation}
U_n z^m = \sum_{I\vDash n} \frac{S^I}{\ell(I)!}\frac{m!}{(m-\ell(I))!}z^{m+n}
        = \sum_{I\vDash n}M_I(m)S^Iz^{m+n}=S_n(mA)z^{m+n}.
\end{equation}

%%%%%%%%%%%%%%%%%%%%%%%%%%%%%%%%%%%%%%%%%%%%%%%%%%%%%%%%%%%%%%%%%%%%%%%%%%%%%%%
%%%%%%%%%%%%%%%%%%%%%%%%%%%%%%%%%%%%%%%%%%%%%%%%%%%%%%%%%%%%%%%%%%%%%%%%%%%%%%%
%%%%%%%%%%%%%%%%%%%%%%%%%%%%%%%%%%%%%%%%%%%%%%%%%%%%%%%%%%%%%%%%%%%%%%%%%%%%%%%
%%%%%%%%%%%%%%%%%%%%%%%%%%%%%%%%%%%%%%%%%%%%%%%%%%%%%%%%%%%%%%%%%%%%%%%%%%%%%%%
%%%%%%%%%%%%%%%%%%%%%%%%%%%%%%%%%%%%%%%%%%%%%%%%%%%%%%%%%%%%%%%%%%%%%%%%%%%%%%%
%%%%%%%%%%%%%%%%%%%%%%%%%%%%%%%%%%%%%%%%%%%%%%%%%%%%%%%%%%%%%%%%%%%%%%%%%%%%%%%
\section{The operad of reduced plane trees}\label{opred}

%%%%%%%%%%%%%%%%%%%%%%%%%%%%%%%%%%%%%%%%%%%%%%%%%%%%%%%%%%%%%%%%%%%%%%%%%%%%%%%
\subsection{A free operad}

We shall now investigate the relation between our Schr\"oder tree expansion
(Section \ref{sec:TreeExp})
and Ecalle's arborification. 
So far, the $S^I$ with $I$  a Schr\"oder pseudocomposition have been
interpreted as elements of the free triduplicial algebra $\SQSym$. They can
also be interpreted as elements of a free operad (see \cite{BH,FC}).  We shall
see that this operad, which is based on reduced plane trees, is also related
to the noncommutative version of the Hopf algebra of formal diffeomorphisms
tangent to identity \cite{BFK}.

The set of reduced plane trees with $n$ leaves will be denoted denoted by
$\text{PT}_n$, and $\text{PT}$ denotes the union
$\bigcup_{n\geqslant 1} \text{PT}_n$. 

The number of leaves of a tree $t$ will be called its degree $d(t)$, and  we
define the grading $\gr(t)$ of a tree as its degree minus 1. In low degrees we
have
\begin{equation}
 \text{PT}_1=\lbrace \arbA \rbrace ,\
 \text{PT}_2=\lbrace \arbB \rbrace ,\ 
 \text{PT}_3=\lbrace \arbCB,\arbCA,\arbCC \rbrace , \dots
\end{equation}

The leaves (in white in the pictures) are also called external vertices whilst the other vertices (in black)
are said to be internal (note that $\arbA$ has no internal vertex). 
For instance, the tree $$t=\arbF$$ has degree $d(t)=8$, grading $\gr(t)=7$ and
$i(t)=4$ internal vertices.

\begin{definition}% (See \cite{FC}) 
The free non-$\Sigma$ operad $\mathcal{S}$ in the category of vector spaces is
the vector space
\begin{equation}
\mathcal{S}=\bigoplus_{n\geqslant 1} \mathcal{S}_n
\ \text{where $\mathcal{S}_n=\C {\rm PT}_n$.}
\end{equation} 
The composition operations 
\begin{equation}
\mathcal{S}_n \otimes \mathcal{S}_{k_1} \otimes \ldots
\otimes \mathcal{S}_{k_n} 
\longrightarrow \mathcal{S}_{k_1 + \ldots + k_n}  
\text{ ($n \geqslant 1$, $k_i \geqslant 1$)}
\end{equation}
map the tensor product of trees $t_0 \otimes t_1 \otimes \ldots \otimes t_n$
to the tree $t_0 \circ (t_1, \ldots, t_n)$ obtained by replacing the leaves of
$t_0$, from left to right, by the trees $t_1, \ldots, t_n$.
\end{definition}

For instance,
\begin{equation}
\arbE \circ \left( \arbA,\arbD,\arbA,\arbB \right)=\arbF
\end{equation}

The tree $\arbA$ of $\text{PT}_1$ is the unit of this composition.
A proof of its associativity can be found in \cite{BH} or \cite{FC}, 
where this operad is called a free $S$-magmatic operad. 
Note that $\mathcal{S}$ is also called the operad of Stasheff polytopes (see
\cite{RH,JLL}) so that the letter $\mathcal{S}$ can stand for Stasheff or
Schr\"oder as well.

%%%%%%%%%%%%%%%%%%%%%%%%%%%%%%%%%%%%%%%%%%%%%%%%%%%%%%%%%%%%%%%%%%%%%%%%%%%%%%%%
\subsection{The group of the operad}

Let $\hat{\mathcal{S}}$ be the completion of the vector space 
$\mathcal{S}$ with respect to the grading $\gr(t)=d(t)-1$.
The group of the operad $\mathcal{S}$ is defined as: 

\begin{definition}
Let 
\begin{equation}
 G_{\rm ncdiff} = \left\{ \arbA + \sum_{n \geqslant 2} p_{(n)} ,
   \hspace{1em} p_{(n)} \in \mathcal{S}_n \right\} \subset \hat{\mathcal{S}}
\end{equation}
endowed with the composition product 
\begin{equation}
 p \circ q = q + \sum_{n \geqslant 2} p_{(n)} \circ \left(
   \underset{n}{\underbrace{q, \ldots, q}} \right) \in G_{{\rm ncdiff}} 
\end{equation}
for $p = \arbA + \sum_{n \geqslant 2} p_{(n)}$ and $q \in G_{{\rm ncdiff}}$. 
\end{definition}

This is indeed a group (see {\it e.g.}, \cite{FC}).
Elements of $G_{{\rm ncdiff}}$  can be described by their coordinates
\begin{equation}
p = \sum_{t \in \text{PT}} p_t t\quad \text{and}\quad q
  = \sum_{t \in \text{PT}} q_t t.
\end{equation}
(with $q_{\arbA} = p_{\arbA} = 1$) so that the coordinates of $r = p \circ q$
are given by
\begin{equation}
r_t = \sum_{ t = t_0 \circ (t_1, \ldots, t_n)} p_{t_0} q_{t_1} \ldots q_{t_n}
\end{equation}
This expression involves the so-called admissible cuts defining the coproduct
in Hopf algebras of the Connes-Kreimer family.
It suggests that the elements of $ G_{{\rm ncdiff}}$ can be interpreted as
characters of the bialgebra defined as follows.

Let
\begin{equation}  
\bar T (\mathcal{S}) = \bigoplus_{p \geqslant 1} \mathcal{S}^{\otimes^p}
\end{equation}
be the reduced tensor algebra over $\mathcal{S}$ 
(whose basis is given by plane forests $f = t_1 \cdot \ldots \cdot t_k$) 
equipped with the coalgebra structure defined on trees by
\begin{equation} 
\tilde{\Delta} (t) = \sum_{
t = t_0 \circ (t_1, \ldots t_n)} t_0 \otimes (t_1 \cdot
   \ldots \cdot t_n) 
\end{equation}
where $\cdot$ means concatenation,
and then extended as an algebra morphism on $\bar T (\mathcal{S})$. 
For example,
\begin{equation}
\tilde{\Delta} \left( \arbCA \right)
 = \arbA \otimes \arbCA + \arbB \otimes \arbB \cdot \arbA
  +\arbCA \otimes \arbA \cdot \arbA \cdot \arbA.
\end{equation}
It is then clear that any $p \in G_{\text{ncdiff}}$ can be identified as the
algebra morphism $\varphi_p$ defined by 
\begin{equation}
\varphi_p (t) = p_t
\end{equation} 
so that, if $r = p \circ q$, then
\begin{equation} 
\varphi_{p \circ q} = \varphi_r = \varphi_p \ast \varphi_q = \mu \circ
   (\varphi_p \otimes \varphi_q) \circ \tilde{\Delta} 
\end{equation}
where $*$ is the usual convolution product for a bialgebra, and $\mu$ is the multiplication of $\C$.
%ICI
Note that if $\bigwedge (t_1 \cdot t_2 \cdot \ldots \cdot t_n)$
(with $n \geqslant 2$) is the tree obtained by grafting the trees
$t_1, \ldots, t_n$ to a common root (in other words,
$\bigwedge (t_1 \cdot t_2 \cdot \ldots \cdot t_n)
 = c_n \circ (t_1, \ldots, t_n)$,
where $c_n$ is the corolla with $n$ leaves) the map $\tilde{\Delta}$ is the
unique algebra map such that
\begin{align}
\widetilde{\Delta} (\arbA) &= \arbA \otimes \arbA,\\
\tilde{\Delta} \left( \bigwedge (t_1, \ldots, t_n) \right)& = \arbA
  \otimes \bigwedge (t_1, \ldots, t_n) + \left( \bigwedge \otimes \text{Id}
  \right) \circ \tilde{\Delta}(t_1 \cdot t_2 \cdot \ldots \cdot
  t_n)
\end{align}

%%%%%%%%%%%%%%%%%%%%%%%%%%%%%%%%%%%%%%%%%%%%%%%%%%%%%%%%%%%%%%%%%%%%%%%%%%%%%%%%
\subsection{The Hopf algebra of reduced plane trees and its characters}

The quotient of the bialgebra $\bar T(\mathcal{S})$ by the relations 
$t \cdot \arbA = \arbA \cdot t = t$ 
is a graded unital algebra $\mathcal{H}_{\text{PT}}$, spanned by  
$\arbA$ and the forests $t_1 \cdot t_2 \cdot \ldots \cdot t_n$
with $t_i \in \cup_{n \geqslant 2} \text{PT}_n$,
with $\arbA$ as unit.

It is a Hopf algebra for the coproduct defined on trees by 
\begin{equation}
\Delta (t) =
(p \otimes p) \circ \tilde{\Delta} (t)
\end{equation}
where $p$ is the projection from $\bar T(\mathcal{S})$ to $\mathcal{H}_{\text{PT}}$. 

In the former example:
\begin{equation}
{\Delta} \left( \arbCA \right)
 = \arbA \otimes \arbCA + \arbB \otimes \arbB  +\arbCA \otimes \arbA.
\end{equation}

It is now clear that the group $G_{\text{ncdiff}}$ is precisely the
group of characters of this Hopf algebra (for $p \in G_{\text{ncdiff}}$,
$\varphi_p (\arbA) = 1$).

This Hopf algebra was first considered in {\cite{MEF}} (with the opposite
coproduct $P\circ\Delta$, where $P(u\otimes v=v\otimes u$)), 
where it is called the Hopf algebra of reduced plane trees
$\mathcal{H}^{\text{red}}_{\text{pl}}$.
The coproduct can be described in terms of admissible cuts of a tree $t \in
\text{PT}$, \emph{i.e.}, (possibly empty) subsets $c$ of edges 
{\it not connected to a leaf\/} with the rule that along any path from the
root of $t$ to any of its leaves, there is at most one edge in $c$. 
The edges in $c$ are naturally ordered from left to right. 
To any admissible cut $c$ corresponds  a unique subforest $P^c(t)$, 
the {\it pruning\/},  concatenation of the subtrees obtained by cutting the
edges in $c$, in the order defined above.
The coproduct can then be defined by:
\begin{equation}
	\Delta(t)=\sum_{c\in\text{Adm }t}R^c(t)\otimes P^c(t),
\end{equation}
where $R^c(t)$ is the trunk, obtained by replacing each subtree of
$P^c(t)$ with a single leaf.

So far, we  have defined  the group of the operad of Stasheff polytopes (or
Schr\"oder trees), and shown  that it coincides with the group of characters
of the Hopf algebra of reduced plane trees.
We shall see that it is also a group of \emph{formal noncommutative
diffeomorphisms} related to the noncommutative Lagrange inversion (see
\cite{NTLag}) and to the noncommutative version of Poincar\'e's equation.

%%%%%%%%%%%%%%%%%%%%%%%%%%%%%%%%%%%%%%%%%%%%%%%%%%%%%%%%%%%%%%%%%%%%%%%%%%%%%%%%
%%%%%%%%%%%%%%%%%%%%%%%%%%%%%%%%%%%%%%%%%%%%%%%%%%%%%%%%%%%%%%%%%%%%%%%%%%%%%%%%
%%%%%%%%%%%%%%%%%%%%%%%%%%%%%%%%%%%%%%%%%%%%%%%%%%%%%%%%%%%%%%%%%%%%%%%%%%%%%%%%
\section{Noncommutative formal diffeomorphisms}\label{ncfd}

%%%%%%%%%%%%%%%%%%%%%%%%%%%%%%%%%%%%%%%%%%%%%%%%%%%%%%%%%%%%%%%%%%%%%%%%%%%%%%%%
\subsection{A group of noncommutative diffeomorphisms}

As pointed out in {\cite{BFK}}, it is possible to consider formal
diffeomorphisms in one variable with coefficients in an associative algebra,
but if this algebra is not commutative, the set of such diffeomorphisms is not
anymore more a group because associativity is broken.
Nevertheless, there is still a noncommutative version of the Fa\`a di Bruno Hopf algebra. 

We can recover a group by regarding the coefficients as well as the variable
as formal noncommutative variables.
Heuristically, let us start with a fixed diffeomorphism $u$ of
\begin{equation} 
G_{\text{diff}} = \{u (z) = z + \sum_{n \geq 1} u_n z^{n + 1} \in
   \mathbb{C}[[z]]\} 
\end{equation}
in the variable $z$ with coefficients $u_n$. 
Consider now that $z$ is replaced by $S_0$ and that each $u_n$ is replaced by
a variable $S_n$.
We get a series 
\begin{equation}
   g_c=S_0 +\sum_{n\geqslant 1} S_n S_0^{n+1} 
\end{equation}
in an infinity of noncommuting variables. 
Nothing prevents us from ``iterating'' $g_c$ as we would do with an ordinary
power series
\begin{equation}
g_c \circ g_c= g_c+\sum_{n\geqslant 1} S_n g_c^{n+1}=S_0 +S_1 S_0^2 +S_2 S_0^3+S_1 (S_1S_0^2)S_0+S_1 S_0 (S_1S_0^2)+...
\end{equation}

Further iterations lead to words in the variables $S_0, S_1, \ldots $ indexed
by Schr\"oder pseudocompositions, which will eventually represent all reduced
plane trees.

Let $S^{\arbA} = S_0$ and, if $t = \bigwedge (t_1, \ldots, t_n)$, $S^t =
S_{n - 1} S^{t_1} \ldots S^{t_n}$. We recover the correspondence with Polish
codes. For example,
\begin{equation}
S^{\arbF}=
S_2S_0S_1S_3S_0S_0S_0S_0S_0S_1S_0S_0=S^{201300000100}
\end{equation}
Identifying trees with their Polish codes, the group $G_{\text{ncdiff}}$ can
be described as
\begin{equation} 
G_{\text{ncdiff}} = \left\{ g = \sum_{t \in \text{PT}} g_t S^t ,
   g_t \in \mathbb{C}, g_{\arbA}=1 \right\} \subset \mathbb{C} \langle \langle S_0, S_1, \ldots \rangle \rangle
\end{equation}
where $\mathbb{C} \langle \langle S_0, S_1, \ldots \rangle \rangle$ 
is the completion of the algebra of polynomials, with respect to the degree in
$S_0$.

If we set
$g = g (S_0, S_1, \ldots) = g (S_0 ; \mathbf{S})$
($\mathbf{S}= S_1, \ldots$), then,

\begin{theorem}
The composition $f \circ g = h$ in $G_{{\rm ncdiff}}$ is given by
\begin{equation}
h (S_0 ; \mathbf{S}) = f (g ; \mathbf{S}). 
\end{equation}
\end{theorem}

In other words we substitute $g$ to the variable $S_0$ in $f$. Graphically we
substitute trees to leaves. 
It is easy to check that this group coincides with the previous one.
If 
\begin{equation}
g = \sum_{t \in \text{PT}} g_t S^t,\  f = \sum_{t \in \text{PT}} f_t
S^t,\ h = f \circ g = \sum_{t \in \text{PT}} h_t S^t \in
G_{{\rm ncdiff}},
\end{equation} 
then $h_t$ is a sum of contributions from $f$ and $g$.
The contributions to $h_t$ can be
\begin{itemize}
\item $f_t^{} S^t$ if we substitute the $S_0$ part of $g$ to any $S_0$ (a
      leaf) of the term $f_t^{} S^t$ of $f$.
  
\item $g_t S^t$ if we substitute the term $g_t S^t$ of $g$ to the $S_0$ part
  of $f$
  
\item $f_{t_0} g_{t_1} \ldots g_{t_n} S^t$ if when substituting (from left
  to right) the terms $g_{t_1} S^{t^1}$, ..., $g_{t_n} S^{t_n}$ to $n$ $S_0$
  variables (leaves) in $f_{t_0} S^{t_0}$, we get the monomial $S^t$.
\end{itemize}
This means that
\begin{equation} 
h_t = \sum_{(t_0 ; t_1 \ldots t_n) = (R^c (t) ; P^c (t))} f_{t_0} g_{t_1}
   \ldots g_{t_n} 
\end{equation}
which is precisely the convolution of characters in $\mathcal{H}_{\text{PT}}$. 

In the sequel, we shall denote by the same letter (for instance $g$) an
element of $G_\text{ncdiff}$ regarded as a series of trees
\begin{equation}
g=\sum_t g_t t,
\end{equation}
as a series of noncommutative monomials
\begin{equation}
g=\sum_t g_t S^t,
\end{equation}
or as the character $g$  sending $t$ on $g_t$.
% (and $f=t_1 \dots t_n$ on $g_{t_1}\dots g_{t_n}$).

%%%%%%%%%%%%%%%%%%%%%%%%%%%%%%%%%%%%%%%%%%%%%%%%%%%%%%%%%%%%%%%%%%%%%%%%%%%%%%%%
\subsection{Inversion in $G_{\rm ncdiff}$ and Lagrange inversion}
\label{sec:invG}
One can compute the compositional inverse of $f_c$ (defined by $f_t=1$ if $t$
is a corolla and $0$ otherwise).
This yields a signed series involving all trees. 

Consider the series 
\begin{equation}
f_c = S^{\arbA} + \sum_{n \geqslant 1} S^{\bigwedge (\arbA^{\cdot n + 1})} . 
\end{equation}
We shall work here in $\bar T(\mathcal{S})$ (where $\arbA$ is
not the unit). Let $i (t)$ be the number of internal
vertices of a tree $t$. The inverse of $f_c$ is then given by 
\begin{equation} 
g_c = \sum_{t \in \text{PT}} (- 1)^{i (t)} S^t 
\end{equation}
since
\begin{align} 
     g_c & =  S_0 +\displaystyle \sum_{n \geqslant 1} 
            \sum_{ t = \bigwedge (t_1 \cdot \ldots \cdot t_{n + 1})\atop t_i \in \text{PT}} 
                    (- 1)^{i (t)} S^{\bigwedge (t_1 \cdot \ldots \cdot t_{n + 1})}\\
     & =  S_0 +\displaystyle  \sum_{n \geqslant 1} 
\sum_{t = \bigwedge (t_1 \cdot \ldots \cdot t_{n + 1})\atop t_i \in \text{PT}} 
(- 1)^{1 + i (t_1) + \ldots i (t_{n + 1})} S_n S^{t_1} \ldots S^{t_{n + 1}}\\
     & =  S_0 - \sum_{n \geqslant 1} S_n g_c^{n + 1}
\end{align}
so that $S_0 = g_c + \sum_{n \geqslant 1} S_n g_c^{n + 1} = f_c \circ g_c$.

In order to establish a link with the noncommutative analogue of Lagrange
inversion (see \cite{NTLag}), we can look for the compositional inverse of
\begin{equation} 
f_L = \left( 1 + \sum_{n \geqslant 1} S_n S_0^n \right)^{- 1} \cdot S_0 \in
   G_{\text{ncdiff}},
\end{equation}
where the exponent $-1$ means here the multiplicative inverse as a formal
power series.
Working with trees, consider the series of trees $L$ such that
\begin{equation}
L = \arbA + \sum_{k \geqslant 1} \bigwedge (L^{\cdot k} \cdot \arbA)
\end{equation}
then $L = \sum L^t t$ with $L^t = 1$ or 0.

The inverse $g_L$ of $f_L$ is $\sum L^t S^t$. Indeed,
\begin{align} 
     g_L & =  S^{\arbA} +\displaystyle  \sum_{k \geqslant 1} S^{\bigwedge (L^{\cdot k} \cdot
     \arbA)}\\
     & =  S_0 +\displaystyle  \sum_{k \geqslant 1} S_k g_L^k S_0\\
     & =  \left(\displaystyle  1 + \sum_{n \geqslant 1} S_n g_L^n \right) S_0^{}
\end{align}
and obviously $f_L \circ g_L = S_0$. 

Apart from $\arbA$, all the trees of $\text{PT}$ occuring in $g_L$ are such
that the rightmost subtree of each internal vertex is a leaf ($S_0$). Let
$\text{PT}_L$ be the set of such trees, and let $\alpha$ be the map sending
$\arbA$ to itself and $t \in \text{PT}_L$ to the tree (with possible unary
internal vertices) obtained by removing all the righmost leaves of its
internal vertices.
We obtain in this way the tree expansion of Section 5.3 in {\cite{NTLag}}.

One can also define $S^t$ for such trees. Now
\begin{equation} 
\alpha (L) = \arbA + \sum_{n \geqslant 1} \bigwedge (\alpha (L)^{\cdot n}) 
\end{equation}
and if we replace trees by their Polish codes, the resulting series $g$
satisfies
\begin{equation} 
g = S_0 + \sum_{n \geqslant 1} S_n g^n, 
\end{equation}
the functional equation considered in \cite{NTLag}. This correspondence will be explained in details  in Section \ref{sec:cat}, using a group morphism from $G_\text{ncdiff}$ to the group $G_\mathcal{C}$ of the Catalan operad.

%%%%%%%%%%%%%%%%%%%%%%%%%%%%%%%%%%%%%%%%%%%%%%%%%%%%%%%%%%%%%%%%%%%%%%%%%%%%%%%%
\subsection{The conjugacy equation}

Let $Y$ be the grading operator on trees ($Y(t)=(d(t)-1)t$), and let 
$q^Y(t)=q^{d(t)-1}t$.
The noncommutative analogue of the conjugacy equation can be written as
\begin{equation}
g(qS_0 ;\mathbf{S})=qg_c(g;\mathbf{S})
\end{equation}
where the initial diffeomorphism is the corolla series. This equation also
reads
\begin{equation}
q^{-1}g(qS_0 ;\mathbf{S})=q^Y g=g_c \circ g.
\end{equation}
It is not difficult to compute the coefficients of the solution
$g=\sum_t c_t(q) t$, noticing that $c_{\arbA}(q)=1$ and, if
$t=\bigwedge(t_1,\dots,t_n)$ then, 
\begin{equation}
q^{d(t)-1} c_t(q)=c_t(q) +c_{t_1}(q)\dots c_{t_n}(q).
\end{equation}
As we have already seen, the coefficients have the closed form
\begin{equation}
c_t(q) = \prod_{v\in t}\frac1{q^{\phi(v)-1}-1}
\end{equation}
where $v$ runs over the internal vertices of $t$ and $\phi(v)$ 
is the number of leaves of the subtree of $v$. 

For example, for
\begin{equation}
t=\arbF,\quad c_t(q)=\frac1{(q^7-1)(q-1)(q^4-1)(q^3-1)}.
\end{equation}

Surprisingly, the same coefficients appear in the commutative case, in
Ecalle's arborified solution of the conjugacy equation, which has been
interpreted in \cite{FM} in terms of characters on the Connes-Kreimer Hopf
algebra $\mathcal{H}_{CK}$ of (non plane) rooted  trees decorated by positive
integers.
We shall see that there is indeed a kind of noncommutative arborification,
which will be eventually explained by a morphism of Hopf algebras.

%%%%%%%%%%%%%%%%%%%%%%%%%%%%%%%%%%%%%%%%%%%%%%%%%%%%%%%%%%%%%%%%%%%%%%%%%%%%%%%%
%%%%%%%%%%%%%%%%%%%%%%%%%%%%%%%%%%%%%%%%%%%%%%%%%%%%%%%%%%%%%%%%%%%%%%%%%%%%%%%%
%%%%%%%%%%%%%%%%%%%%%%%%%%%%%%%%%%%%%%%%%%%%%%%%%%%%%%%%%%%%%%%%%%%%%%%%%%%%%%%%
\section{Commutative versus noncommutative}\label{cvsnc}

%%%%%%%%%%%%%%%%%%%%%%%%%%%%%%%%%%%%%%%%%%%%%%%%%%%%%%%%%%%%%%%%%%%%%%%%%%%%%%%%
\subsection{Commutative diffeomorphisms and the Connes-Kreimer algebra}

In this section we recall briefly how certain (commutative) formal
diffeomorphisms can be obtained as characters of a Hopf algebra (see \cite{FM}
and \cite{smf}), in particular the solution of the conjugacy equation
\begin{equation}
qu(h(z))=h(qz)
\end{equation}
where $u$ and $h$ are formal diffeomorphisms of $G_{\text{diff}}$ (tangent to
the identity). 
The use of trees to encode diffeomorphisms appears in \cite{JESNAG} and is
related to differential operators indexed by trees, an idea originally due to
Cayley.

We shall rely upon the references {\cite{ck0}}, {\cite{fo1}} and {\cite{fo2}},
except that we use the opposite coproduct, in order to avoid antimorphisms. 

A rooted tree\footnote{Not to be confused with rooted plane trees, of which Schr\"oder trees are a special case.} 
$T$ is a connected and simply connected set of oriented edges
and vertices such that there is precisely one distinguished vertex (the root)
with no incoming edge.
A forest $F$ is a (commutative) monomial in rooted trees. 

Let $l(F)$ be the number of vertices in $F$. One can decorate a forest
by $\mathbb{N}^{\ast}$, that is, with each vertex $v$ of $F$, we
associate an element $n_v$ of $\mathbb{N}^{\ast}$. 
We denote by
$\mathcal{T}_{\mathbb{N}}$ (resp. $\mathcal{F}_{\mathbb{N}}$) the set of
decorated trees (resp. forests). It includes the empty tree, denoted by
$\emptyset$.
As for sequences, if a forest $F$ is decorated by $n_1, \ldots,
n_s$ ($l (F) = s$), we write
\begin{equation} 
|F|= n_1 + \ldots + n_s.
\end{equation}

For $n$ in $\mathbb{N}^{\ast}$, the operator $B_n^+$ associates with a forest
of decorated trees the tree with root decorated by $n$ connected to the roots
of the forest : $B_n^+ (\emptyset)$ is the tree with one vertex decorated by
$n$. For example, 
\begin{equation}
 B_n^+ \left(
\begin{arb}
\rd{n_4} child{\va{n_5}};
\end{arb} 
\begin{arb}
\rd{n_1} child{\vl{n_2}} child{\vr{n_3}};
\end{arb}
\right)=
\begin{arb}
\rd{n} child{\vl{n_1} child{\vl{n_2}} child{\va{n_3}}}
       child{\vr{n_4} child{\vr{n_5}}};
\end{arb} 
\end{equation}

The linear span $\mathcal{H}_{\text{CK}}$ of $\mathcal{F}_{\mathbb{N}}$ is
the graded Connes-Kreimer Hopf algebra of trees decorated by
$\mathbb{N}^{\ast}$ for the product
\begin{equation} 
\pi (F_1 \otimes F_2) = F_1 F_2 
\end{equation}
and the unit $\emptyset$. 

The coproduct $\Delta$  can be defined by
induction 
\begin{align}
 \Delta (\emptyset)& = \emptyset \otimes \emptyset,\\ 
\Delta (T_1 \ldots T_k)& = \Delta (T_1) \ldots \Delta (T_k),\\ 
\Delta (B_n^+ (F))& =
  \emptyset  \otimes B_n^+ (F) + ( B_n^+ \otimes \text{Id}) \circ \Delta (F).
\end{align}

There exists a combinatorial description of this coproduct (see {\cite{fo1}}).
For a given tree $T \in \mathcal{T}_{\mathbb{N}}$, an admissible cut $c$ is a
subset of its vertices such that, on the path from the root to an element of
$c$, no other vertex of $c$ is encountered. For such an admissible cut, $P^c
(T)$ is the  product of the subtrees of $T$ whose roots are in $c$ and $R^c
(T)$ is the remaining tree, once these subtrees have been removed. With these
definitions, for any tree $T$, we have
\begin{equation} 
\Delta (T) = \sum_{c \hspace{1em} \text{adm} .} R^c (T) \otimes P^c (T) .
\end{equation}
For example,
\begin{equation*}
 \Delta \left( \begin{arb} \rd{n_1} child{\vl{n_2}} child{\vr{n_3}}; \end{arb}
        \right)
 =
\begin{arb} \rd{n_1} child{\vl{n_2}} child{\vr{n_3}}; \end{arb}
 \otimes \emptyset
 + \begin{arb} \rd{n_1}  child{\vr{n_3}}; \end{arb}
 \otimes \bullet_{n_2}
 + \begin{arb} \rd{n_1} child{\vl{n_2}} ; \end{arb}
 \otimes
 \bullet_{n_3} + \bullet_{n1}\otimes \bullet_{n_2}\bullet_{n_3} 
 + \emptyset
 \otimes \begin{arb} \rd{n_1} child{\vl{n_2}} child{\vr{n_3}}; \end{arb}
\end{equation*}

%For a forest $F$ in $\mathcal{F}_{\mathbb{N}}$ we remind that the symmetry %factor (\cite{fo2}) of $F$
%is defined by :
%\begin{enumerate}
%  \item $s ((\eta)) = 1$ ;
%  
%  \item $s (B^+_n (F)) = s (F)$ ;
%  
%  \item $s (T_1^{a_1} \ldots T_k^{a_k}) = s (T_1)^{a_1} \ldots s (T_k)^{a_k}
%  a_1 ! \ldots a_k$! if $T_1$,...,$T_k$ are distinct rooted trees.
%\end{enumerate}

\begin{definition} \label{defAT} 
Given  a formal diffeomorphism 
\begin{equation}
u(z)=z+\sum_{n\geqslant 1}u_n z^{n+1},
\end{equation}
we associate with any tree $T$ (see \cite{FM,smf,Men06,Men07})
a monomial $A_T(z)$ recursively defined as follows:
\begin{itemize}
    \item For the empty tree $A_{\emptyset}(z)=z$,
    \item If $T = B_+^n (\emptyset)$ then $A_T(z)=u_nz^{n+1}$, 
    \item If $T = B_+^n (F)$ where $F = T_1^{a_1} \ldots T_k^{a_k}$ is a non
empty product  of $k$ distinct decorated trees, with multiplicities $a_1,
\ldots, a_k$ ($a_1 + \ldots + a_k = s$), then
\begin{equation}
    A_T(z) =  \frac{1}{a_1 ! \ldots a_k !}
       A_{T_1}^{a_1}(z) \ldots A_{T_k}^{a_k}(z) 
       \left(\partial_z^{a_1+\dots+a_k} u_n z^{n+1}\right)
\end{equation}
\end{itemize}
\end{definition}

The main result is then:
\begin{theorem}
The map $\rho_u$ sending a character $\varphi$ of
$\mathcal{H}_{\text{CK}}$ to the formal diffeomorphism of $G_{{\rm diff}}$ 
\begin{equation}
\rho_u(\varphi)(z)=\sum_{T\in \mathcal{T}_{\mathbb{N}}}\varphi(T)A_T(z)
\end{equation}
is a group homomorphism from the group of characters $G_{CK}$ of $\mathcal{H}_{\text{CK}}$
to $G_{{\rm diff}}$.
\end{theorem}

See \cite{FM,smf,Men06,Men07} for proofs. 

\bigskip 
Going back to the equation $u(h(z))=q^{-1}h(qz)$, one can observe that
\begin{itemize}
\item $u=\rho_u(\varphi_0)$, where $\varphi_0$ is the character 
      given by $\varphi_0(T)=1$ (resp. $0$)
      if $T=\emptyset$ or $T=B_+^n (\emptyset)$ (resp. otherwise).
\item If the conjugating $h$ is given by a character $\theta$ 
($h=\rho_u(\theta)$) then $q^{-1}h(qz)$ is given by the character
 $\theta \circ q^Y$ where $q^Y(F)=q^{|F|} F$.
\end{itemize}
 
Therefore, the conjugacy equation can be lifted to the character equation
\begin{equation}
\varphi_0\ast\theta=\theta \circ q^Y.
\end{equation} 
This equation is easily solved.  
For a tree $T = B_+^n (T_1 \ldots T_s)$, we get
\begin{equation}
(q^{|T|} - 1) \theta(T) = \theta(T_1) \ldots \theta(T_s)
\end{equation}
so that 
\begin{equation}
\theta(T)=\prod_{v\in T}\frac1{q^{|T_v|}-1}
\end{equation}
where $T_v$ is the subtree of $T$ whose root is the vertex $v$. 

A more explicit expression can be found in \cite{smf}. 
Such ``arborified'' expressions are useful for analysis since they allow to
prove the analyticity of the conjugating map under some diophantine conditions
on $q$ (ensuring a geometric growth of the numbers $\theta(T)$, see
{\it e.g.}, \cite{smf}).

%%%%%%%%%%%%%%%%%%%%%%%%%%%%%%%%%%%%%%%%%%%%%%%%%%%%%%%%%%%%%%%%%%%%%%%%%%%%%%%%
\subsection{Relating $\mathcal{H}_{\text{CK}}$ and $\mathcal{H}_{\text{PT}}$}

The similarity of the coefficients $c_t(q)$ and $\theta(T)$ suggests a link 
between both versions of the conjugacy equation,
which turns out to be understandable at the level of Hopf algebras.

\begin{theorem}
Let $\sk$ (for ``skeleton'') be the map defined from $PT$ to
 $\mathcal{T}_{\mathbb{N}}$ by $\sk(\arbA)=\emptyset$ and, 
if $t = \bigwedge (t_1, \ldots, t_n)$ ($n\geqslant 2$), 
then $\sk(t)=B_+^{n-1} (\sk(t_1)\ldots \sk(t_n))$. 

This map extends naturally to an algebra morphism from 
$\mathcal{H}_{\text{PT}}$ to $\mathcal{H}_{\text{CK}}$ which is actually a
Hopf algebra morphism.
\end{theorem}

For example if $t=\arbF$,
\begin{equation}
T=\sk \left( t \right)=
\begin{arb}
\rd{2} child{\vl{1}} child{\vr{1} child{ \va{3}}};
\end{arb}
\end{equation}

The proof follows immediately by comparing the recursive (or combinatorial)
definition of the coproducts.

This Hopf morphism induces a group morphism $\sk^{\ast}$ sending a character
$\varphi\in G_{\text{CK}}$ 
to the character $\sk^{\ast}(\varphi)=\varphi\circ \sk$. 
It sends the character $\theta$ to the character defined by
\begin{equation}
g=\sum_t c_t(q)t.
\end{equation}
In other words, for any reduced plane tree $t$, $c_t(q)=\theta(\sk(t))$. 

%%%%%%%%%%%%%%%%%%%%%%%%%%%%%%%%%%%%%%%%%%%%%%%%%%%%%%%%%%%%%%%%%%%%%%%%%%%%%%%%
\subsection{The final diagram}

Let us summarize the situation:
\begin{itemize}
\item We have a noncommutative analogue of the conjugacy equation whose
solution is a character on $\mathcal{H}_{\text{PT}}$ defined by the
coefficients $c_t(q)$,
\item In the commutative case, the solution can be computed as a character on
$\mathcal{H}_{\text{CK}}$,
\item Both characters are related by the group morphism $\sk^{\ast}$.
\end{itemize}

There is also a morphism from $G_{\text{ncdiff}}$ to $G_{\text{diff}}$,
so that we also get the solution in the commutative case. 
For $u \in G_{\text{diff}}$, the algebra map $\alpha_{u}$ sending $S_0$ to $z$
and $S_n$ to $u_n$ defines a group morphism and, if we denote by
$g_c \in G_{\text{ncdiff}}$ the sum of all corollas, then
$\alpha_{u} (g_c) = u$. 
This morphism also sends the solution of the noncommutative conjugacy equation
to the solution of the commutative one.

The final picture is then:
\begin{theorem}
The diagram
\begin{equation}
\begin{array}{ccc}
G_{\text{CK}} & \overset{\text{\sk}^{\ast}}{\longrightarrow}
 & G_{\text{ncdiff}}\\
  &\underset{\rho_{u}}{\searrow} & \downarrow \text{}_{\alpha_{u}} \\
     & & G_{\text{diff}}
\end{array}
\end{equation}
is commutative.
\end{theorem}

\Proof
Consider a character $\varphi \in G_{\text{CK}}$. Then, on the one hand
\begin{equation} 
\rho_{u} (\varphi) = \sum_T \varphi(T)A_T(z)
\end{equation}
and, on the other hand, as 
$\phi = \sk^{\ast}(\varphi)$ ($\phi (t) =  \varphi(\sk (t)$), 
\begin{equation} 
\alpha_u  (\phi) = \alpha_u \left( \sum_t \phi(t) S^t \right)
 = \sum_T \varphi(T) \alpha_u \left( \sum_{\sk (t) = T} S^t \right). 
\end{equation}
The above diagram commutes if, for any tree $T$,
\begin{equation} 
A_T(z) = \alpha_u \left( \sum_{\sk(t) = T} S^t \right). 
\end{equation}
The result is obvious for $T=\emptyset$, and for any tree
$C_n=B^+_n(\emptyset)$, since $A_{C_n}(z)=u_n z^{n+1}$ and the only reduced
tree with such a skeleton is the corolla $c_{n+1}$:
\begin{equation}
\alpha_u (S^{c_{n+1}})=\alpha_u (S_nS_0^{n+1})=u_n z^{n+1}=A_{C_n}(z).
\end{equation}

Now, from Definition~\ref{defAT}, if $T = B^+_n (T_1^{a_1} \ldots T_k^{a_k})$
where $T_1$,...,$T_k$ are distinct rooted trees, then
\begin{equation}
 A_T (z) = \left(\begin{array}{c}
     n + 1\\
     a_1, a_2, \ldots, a_k
   \end{array}\right) A_{T_1}^{a_1} (z) \ldots A_{T_k}^{a_k} (z) u_n
   z^{n + 1 - a_1 - \ldots - a_k} .
\end{equation}
Let 
\begin{equation}
\text{$X_0 = S_0$ and}\ X_i =  \sum_{\sk (t) = T_i} S^t. 
\end{equation}
If $W (a_0, a_1, \dots, a_k)$ is the set of all possible words obtained by
concatenating $a_i$ copies of $X_i$, then
\begin{equation}
\sum_{\text{\sk} (t) = T} S^t  = \sum_{w \in W (n + 1 - a_1 - \ldots - a_k,
   a_1, \ldots, a_k)} S_{n}w 
\end{equation}
since a tree has skeleton $T$ if and only if it can be written
$c_n\circ(t_1,\dots t_n)$ for the operadic composition, where the $n$-tuple
$(t_1,\dots t_n)$ contains exactly $a_i$ trees with skeleton $T_i$ and
$n+1-a_1 -\dots -a_k$ trees $\arbA$.
There are exactly 
$\left(\begin{array}{c}
     n + 1\\
     a_1, a_2, \ldots, a_k
   \end{array}\right)$ such $n$-tuples, and, by induction,

\begin{align}
\alpha_u \left( \sum_{\sk (t) = T} S^t \right)& =A_{T_1}^{a_1}(z)\ldots
   A_{T_k}^{a_k} (z) u_n z^{n + 1 - a_1  \ldots - a_k} \left( \sum_{w
   \in W (n + 1 - a_1  \dots - a_k, a_1, \ldots, a_k)} 1 \right)\\
 & = A_T (z)\nonumber
\end{align}
\qed

%%%%%%%%%%%%%%%%%%%%%%%%%%%%%%%%%%%%%%%%%%%%%%%%%%%%%%%%%%%%%%%%%%%%%%%%%%%%%%%%

\section{The Catalan operad}\label{sec:cat}
\def\Cat{{\mathcal C}}
\def\TT{{\mathcal T}}
The preceding considerations can be repeated almost word for word for
the free operad on one generator in each degree
\begin{equation}
\Cat = \bigoplus_{n\ge 1}\C \TT_n
\end{equation}
where $\TT_n$ denotes the set of all plane trees on $n$ vertices:
\begin{equation}
 \TT_1=\lbrace \carbA \rbrace ,\ 
 \TT_2=\lbrace \carbUn \rbrace ,\ 
 \TT_3=\lbrace \carbBA, \carbBB \rbrace ,\ 
 \TT_4=\lbrace \carbCA,\carbCB,\carbCC,\carbCD, \carbCE \rbrace , \dots
\end{equation}
endowed with the same composition as in $\mathcal{S}$.
We shall call it the {\it Catalan operad}, although
it probably has other names \cite{RH,FC}.

General plane trees can be represented by their Polish codes, as monomials
$S^I$, where $S_n$ is now the symbol of an $n$-ary operation, that is $S^{\carbA}=S_0$ and $S^{B_{+}(t_1,\dots,t_k)}=S_k S^{t_1}\dots S^{t_k}$. For instance,
\begin{equation}
S^{\carbF}=S^{3100200}.
\end{equation}
  As observed
in \cite{NTLag}, these $I$ are also the evaluation vectors on nondecreasing
parking functions.

The discussion of Lagrange inversion and related problems given in \cite{NTLag}
can be interpreted as calculations in the group $G_\Cat$ of this operad. 
The functional equation
\begin{equation}\label{eq:lagtree}
g = \sum_{n\ge 0}S_ng^n
\end{equation}
can be rewritten as
\begin{equation}
(S_0-S^{10}-S^{200}-\cdots)\circ g = S_0
\end{equation}
so that it amounts to inverting the series
\begin{equation}
f = S_0-\sum_{n\ge 1}S^{ n 0^n}
\end{equation}
in the group $G_\Cat$. 
The  relation with the computations in Section \ref{sec:invG} can be elucidated with the help of a surjective group morphism 
from $G_{\rm ncdiff}$ to $G_\Cat$. 
Recall that $\text{PT}_L$ is the set of Schr\"oder trees such that the rightmost subtree of each internal vertex is a leaf. 
Let $\alpha$ be the map which sends $\arbA$ to $\carbA$ and any other tree $t$ of $\text{PT}$ :
\begin{itemize}
\item to the the plane tree obtained by removing all the rightmost leaves of its
internal vertices if $t$ is in $\text{PT}_L$,
\item to 0 otherwise.
\end{itemize} 

This map is surjective on plane trees and induces a linear map on monomials in $S_0, S_1, ...$ 
which happens to be a group morphism from $G_{\rm ncdiff}$ to $G_\Cat$. If we still denote  by $\alpha$ this morphism, then, since
\begin{align}
\alpha( g_c) = \alpha \left(S^{\arbA} + \sum_{n \geqslant 1} S^{\bigwedge (\arbA^{\cdot n + 1})} \right) &=  S_0+\sum_{n\ge 1 }S^{ n 0^n}\in G_\Cat, \\
\alpha(g_L) = \alpha \left( \left( 1 + \sum_{n \geqslant 1} S_n S_0^n \right)^{- 1} S_0\right) &= 
S_0-\sum_{n\ge 1}S^{ n 0^n}\in G_\Cat.
\end{align}
we can obtain the composition inverse of these series of $G_\Cat$ as $\alpha(f_c)$ and $\alpha(f_L)$,  which gives back the formulas of  \cite{NTLag}.

The functional equation for $g$  can also be written
\begin{equation}\label{eq:lagbin}
g = S_0+\Omega g\cdot g =: S_0 + B(g,g)
\end{equation}
Each plane tree in the solution of \eqref{eq:lagtree} corresponds to a unique binary tree $B_T(S_0)$
in the solution of \eqref{eq:lagbin}.
This induces a bijection between plane trees and binary trees: writing (see (\ref{def:omeg}))
\begin{equation}
B(S^I,S^J) = \Omega S^IS^J,
\end{equation}
there is a unique way to decompose a plane tree $S^I$ on $n$ vertices as
\begin{equation}
S^I=B(S^{I_1},S^{I_2}) 
\end{equation}
so that recursively
\begin{equation}
S^I=B_T(S_0,\ldots,S_0)
\end{equation}
for a unique binary tree $T$ with $n-1$ internal vertices. For example,
\begin{align}
S^{\carbUn} = S^{10} &= B(S_0,S_0),\\
S^{\carbBA} = S^{200} &= B(S^{10},S_0)=B(B(S_0,S_0),S_0),\\
S^{\carbBB} = S^{110}&= B(S_0,S^{10}) = B(S_0,B(S_0,S_0))
\end{align}
We recover in this way the classical rotation correspondence.

In fact, if, in the one to one correspondence with plane tree, $S^{I_1}=S^{t_1}$ and $S^{I_2}=S^{t_2}$, the tree corresponding to $S^I=B(S^{I_1},S^{I_2})=\Omega S^IS^J$ is $t=B_+(B_-(t_1),t_2)$. Using this trick we get for instance
\begin{equation}
S^{3100200}=S^{\carbF}=B(S^{\carbCB},S^{\carbBA})=B(B(S^{\carbBB},S^{\carbA}),B(S^{\carbUn}, S^{\carbA}))
\end{equation}
that corresponds finally to the binary tree $B(B(B(S_0,B(S_0,S_0)),S_0),B(B(S_0,S_0),S_0))$:
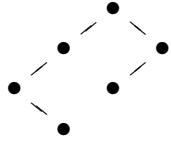
\def\carbX#1#2#3#4#5#6{
\xymatrix@R=0.1cm@C=2mm{
 & & & {\GrTeXBox{#1}}\arx1[ld]\arx1[rd]\\
% ligne 2
 & & {\GrTeXBox{#2}}\arx1[ld] & & {\GrTeXBox{#3}}\arx1[ld] \\
% ligne 3
 & {\GrTeXBox{#4}}\arx1[rd] && {\GrTeXBox{#5}} \\
 & & {\GrTeXBox{#6}}
}
}
\begin{equation}
\carbX{\bullet}{\bullet}{\bullet}{\bullet}{\bullet}{\bullet}
\end{equation}

Another question which can be investigated in this context is the formal solution
of the generic differential equation
\begin{equation}
\frac{dx}{ds} = f(x(s)),\quad x(0)=x_0.
\end{equation}
Rather than stating Cayley's formula for $x^{(k)}$ in terms of rooted trees and
derivatives, we shall write down a noncommutative version involving plane trees and
the coefficients of the generic power series $f$. Assuming without loss of generality
that $x_0=0$, we can look for a series $X(s)\in G_\Cat$ satisfying
\begin{equation}\label{eq:ncdiff}
\frac{dX}{ds} = \sum_{n\ge 0} S_n X(s)^n.
\end{equation}
Thus,
\begin{equation}
X(s) = S_0s+ \sum_{n\ge 1} S_n \int_0^s X(u)^n du =: \sum_{n\ge 1} X_n \frac{s^n}{n!}
\end{equation}
and solving iteratively as usual, we get successively
\begin{align}
X_1 &= S_0,\\
X_2 &= S^{10},\\
X_3 &= 2S^{200}+S^{110},\\
X_4 &= 6S^{3000}+3S^{2100} + 3S^{2010} + 2S^{1200} +S^{1110}.
\end{align}
Identifying as before trees and their Polish codes, 
\begin{equation}
X_n = \sum_{t\in \TT_n}c_t S^t
\end{equation}
and setting
\begin{equation}
F_n(Y_1\ldots,Y_n) = S_n\int_0^sY_1(u)\cdots Y_n(u)du
\end{equation}
we have
\begin{equation}
X_{n+1} = \sum_{k=1}^n \sum_{i_1+\cdots+i_k=n}F_k\left(X_{i_1}\frac{s^{i_1}}{i_1!},\ldots,X_{i_k}\frac{s^{i_k}}{i_k!}\right),
\end{equation}
which gives for the coefficient $c_t$ of $t=B_+(t_1,\ldots,t_k)\in\TT_{n+1}$
\begin{equation}
c_t\frac{s^{n+1}}{(n+1)!}= \int_0^s c_{t_1}\cdots c_{t_k}\frac{u^{|t_1|+\cdots+|t_k|}}{|t_1|!\cdots |t_k|!}du
\end{equation}
so that
\begin{equation}
c_t = {n\choose |t_1|,\ldots,|t_k|}c_{t_1}\cdots c_{t_k}
\end{equation}
which is clearly the recursion for the number of decreasing (or increasing) labellings
of $t$, also given by the hook-length formula
\begin{equation}
c_t = (n+1)!\prod_{v\in t}\frac1{h_v}
\end{equation}
where $h_v$ is the number of nodes of the subtree with root $v$. For instance, for the tree $\carbCB$ that corresponds to the monomial $S^{2100}$, there are 3 decreasing labelings: 
\begin{equation}
\begin{arb}
\rd{4} child{\vl{3}} child{\vr{2} child{ \va{1}}};
\end{arb},\begin{arb}
\rd{4} child{\vl{2}} child{\vr{3} child{ \va{1}}};
\end{arb}, \begin{arb}
\rd{4} child{\vl{1}} child{\vr{3} child{ \va{2}}};
\end{arb}.
\end{equation}

Replacing each $S_k$ by $\frac{f^{(k)}}{k!}$, we recover Cayley's formula for
$x^{(n)}$, 
\begin{equation}
 {d^n x\over ds^n} = \sum_{|t|=n} a(t) \delta_t,
\end{equation}
where
\begin{equation}
a(t)= {|t|!\over t! |S_t|},
\end{equation}
$S_t$ being the symmetry group of $t$, 
\begin{equation}
 B_+(t_1,\dots,t_n)! = |B_+(t_1,\dots,t_n)| \cdot t_1! \cdots t_n!,\quad  \bullet ! =1.
\end{equation}
and the elementary differentials are defined by \cite{Bu}
\begin{equation}
 \delta_\bullet^i= f^i, \,\,\, \delta^i_{B_+(t_1,\dots,t_n)} = \sum_{j_1,\dots,j_n=1}^N (\delta^{j_1}_{t_1} \cdots \delta^{j_n}_{t_n})\partial_{j_1} \cdots \partial_{j_n} f^i
\end{equation}
In particular, the solution is given by
\begin{equation}
\displaystyle x(s) = x_0 + \sum_{t}  {s^{|t|}\over |t|!} a(t) \delta_t(0)
\end{equation}
For example,
\begin{equation}
 x^{(4)} = f^{\prime\prime\prime}(f,f,f)  + 3f^{\prime\prime}(f,f^\prime(f))   +   f^\prime(f^{\prime\prime}(f,f)) +f^\prime(f^\prime(f^\prime(f))).
\end{equation}
Note that with the interpretation of $S_n$ as an $n$-linear operation, our formal calculations are valid for $x\in\R^N$: we can write
the Taylor expansion of $f$ as
\begin{equation}
f(x) = F_0 + F_1(x)+F_2(x,x)+F_3(x,x,x) + \cdots
\end{equation}
without expliciting the expression of $F_n$ in terms of partial derivatives.

Once again, the functional equation \eqref{eq:ncdiff} can be recast as a quadratic fixed point
problem:
\begin{equation}
\frac{dX}{ds} = S_0 + (S_1 + S_2X(s) +S_3X^2(s)+\cdots)X(s) = S_0 + (\Omega X'(s))X(s)
\end{equation}
so that
\begin{equation}
X(s) = S_0s + \int_0^s \Omega X'(u)\cdot X(u)du = S_0s + \B(X(s),X(s))
\end{equation}
The bilinear map $\B$ acts on trees by
\begin{equation}
\B\left (S^I \frac{s^i}{i!}, S^J\frac{s^j}{j!}\right) ={i+j\choose i,j}\Omega S^{IJ} \frac{s^{i+j+1}}{(i+j+1)!}
\end{equation}

%%%%%%%%%%%%%%%%%%%%%%%%%%%%%%%%%%%%%%%%%%%%%%%%%%%%%%%%%%%%%%%%%%%%%%%%%%%%%%%%

\section{Appendix: numerical examples}\label{sec:ex}

In the case of Lagrange inversion, comparison between the formal 
noncommutative solution and numerical examples 
(specializations of the alphabet, or characters) leads to interesting insights.
We shall give here a (short) list of known workable examples.

\subsection{Warmup: $A=1$}

The alphabet $A=1$ is defined by $S_n(1)=1$ for all $n$.
In this case,
\begin{equation}
\phi(z)= \frac{qz}{1-z} = qz\sigma_z(1)
\end{equation}
is a M\"obius transformation, and it is trivial to conjugate it
to its linear part when $q\not=1$. However, it is a good exercise
to work out the series solution. We have
\begin{equation}
S_n(m)=\binom{n+m-1}{n}
 \ \text{so that}\ 
 L(S^I)(1)=\frac{n!}{i_1!\cdots i_r!}=n!S^I(\E)
\end{equation}
where $\E$ is defined by $S_n(\E)=1/n!$. Hence,
\begin{equation}
f(z)= \int_0^\infty
      e^{-t} L\left( z\sigma_{zt}\left(\frac{\E}{1-q}\right)\right) dt
 = \frac{z}{1-\frac{z}{1-q}}.
\end{equation}

\subsection{The logistic map: $A=-1$}

The logistic map is defined by
\begin{equation}
\phi(z) = qz(1-z) = qz\sigma_z(-1).
\end{equation}
Indeed, by definition, $\sigma_z(-1)$ is the inverse
of $\sigma_z(1)$, so that $S_1(-1)=-1$ and $S_n(-1)=0$
for $n>1$. 

The recurrence for %$h_n=
$ g_n(-1)$ is here
\begin{equation}
(1-q^n)g_n = \sum_{k=0}^{n-1}g_kg_{n-1-k}
\end{equation}
yielding
\begin{equation}
g_1=\frac1{1-q},\quad g_2=\frac2{(1-q)(1-q^2)},\quad g_3=\frac{5+q}{(1-q)(1-q^2)(1-q^3)},\ldots
\end{equation}
the numerator being a $q$-analogue of $n!$.

For $q=-2,2,4$, these series have explicit forms in terms of elementary functions:
\begin{align}
q=-2:\ & f(z)     =\frac{\sqrt{3}}{6}\left(  2\pi-3\arccos\left( z-\frac{1}{2}\right)  \right),
\quad h(z) =\frac{1}{2}-\cos\left(  \frac{2z}{\sqrt{3}}+\frac{\pi}{3}\right) ,\nonumber\\
q=2:\ & f(z)= -\frac{1}{2}\ln\left(  1-2z\right),\quad h(z)  =\frac{1}{2}\left(  1-e^{-2z}\right) ,\nonumber\\
q=4:\ & f(z)= \left(  \arcsin\sqrt{z}\right)^{2},\quad h(z) =\left(  \sin\sqrt{z}\right)^{2} .
\label{eq:exact}
\end{align}
Numerical investigations, including a conjecture for the radius of convergence in the general case,
can be found in \cite{CZ1,CZ2}.

\subsection{The Ricker map: $A=\E$}

This case corresponds to
\begin{equation}
\phi(z)=qze^z.
\end{equation}
No closed expression is known for $f$ or $g$, but a numerical study can be found in \cite{CZ1}.

We have
\begin{equation}
S_n(m\E)=\frac{m^n}{n!},\ \text{so that}\ L(S^I)=\frac{1^{i_1}(i_1+1)^{i_2}(i_1+i_2+1)^{i_3}\cdots(i_1+\cdots i_{r-1}+1)^{i_r}}{i_1!i_2!\cdots i_r!}
\end{equation}
and we can compute
\begin{equation}
f_1=\frac1{1-q},\quad f_2=\frac{3q+1}{2!(1-q)(1-q^2)},\quad f_3=\frac{16q^3+11q^2+8q+1}{3!(1-q)(1-q^2)(1-q^3)},\ldots
\end{equation}
The numerators are $q$-analogues of $(n!)^2$, whose combinatorial interpretation requires further
investigations.

%%%%%%%%%%%%%%%%%%%%%%%%%%%%%%%%%%%%%%%%%%%%%%%%%%%%%%%%%%%%%%%%%%%%%%%%%%%%%%%
%%%%%%%%%%%%%%%%%%%%%%%%%%%%%%%%%%%%%%%%%%%%%%%%%%%%%%%%%%%%%%%%%%%%%%%%%%%%%%%
%%%%%%%%%%%%%%%%%%%%%%%%%%%%%%%%%%%%%%%%%%%%%%%%%%%%%%%%%%%%%%%%%%%%%%%%%%%%%%%

\end{document}